\documentclass{amsart}
\usepackage{amssymb}
\usepackage{enumitem}
\usepackage[T1]{fontenc}
\usepackage{lmodern}
\usepackage{microtype}
\usepackage[margin=1in]{geometry}
\usepackage[colorlinks=true,linkcolor=black,citecolor=black,urlcolor=blue]{hyperref}
\usepackage{silence}
\theoremstyle{definition}
\newtheorem{thm}{Theorem}[section]
\newtheorem{prop}[thm]{Proposition}
\newtheorem{lem}[thm]{Lemma}
\newtheorem{cor}[thm]{Corollary}

\setlist[enumerate]{label=(\alph*),leftmargin=1.8em,itemsep=3pt}

\newcommand{\mr}{\mathrm}
\newcommand{\bb}{\mathbb}
\newcommand{\mb}{\mathbf}
\newcommand{\mc}{\mathcal}
\newcommand{\bE}{\bb E}
\newcommand{\bP}{\bb P}
\newcommand{\bR}{\bb R}
\newcommand{\bZ}{\bb Z}
\newcommand{\Fhat}{\widehat F}
\newcommand{\Sur}{\mr{Sur}}
\newcommand{\Hom}{\mr{Hom}}
\newcommand{\Aut}{\mr{Aut}}
\newcommand{\Gen}{\mr{Gen}}
\newcommand{\supp}{\mr{supp}}
\newcommand{\Pois}{\mr{Pois}}
\newcommand{\Bin}{\mr{Bin}}
\newcommand{\Tr}{\mr{Tr}}
\newcommand{\Var}{\mr{Var}}
\newcommand{\op}{\mr{op}}
\newcommand{\dist}{\xrightarrow{\mr{dist}}}

\title[Sharp universality thresholds for random profinite groups]
{Sharp universality thresholds for random profinite groups}
\author{Jiwan Jung}
\address{Department of Mathematics, Pohang University
of Science and Technology, Pohang 37673,
Republic of Korea}
\email{guinipig123@postech.ac.kr}

\author{Jungin Lee}
\address{Department of Mathematics, Ajou University,
Suwon 16499, Republic of Korea}
\email{jileemath@ajou.ac.kr}

\begin{document}

\begin{abstract}
We determine sharp universality thresholds for quotients of the free profinite group $\widehat F_n$ on $n$ generators by $n+u$ independent random relations, where $u$ is a fixed nonnegative integer. For products of mutually independent uniformly coset-balanced factors, convergence to the Liu--Wood measure $\mu_u$ holds when the minimum number of factors grows faster than $\log n$, which is a uniformly sharp condition. For random words of fixed length $\ell_n$, convergence holds exactly when $\ell_n$ is eventually odd and $\ell_n-\log n\to\infty$. 
We also prove universality for randomly ordered products of independent integer powers under the balance condition $n\alpha_{n,p}/\log n\to\infty$ at each fixed prime, recovering finite-prime cokernel laws by abelianization.
\end{abstract}

\maketitle
\vspace{-5mm}

\section{Introduction}\label{Sec1}

Let $\Fhat_n$ be the free profinite group on $n$ generators $x_1,\ldots,x_n$, and let $u\ge0$ be an integer. For independent random relations $r_{n,1},\ldots,r_{n,n+u}\in\Fhat_n$, set
$$
\mc G_{n,u}:=\Fhat_n/\langle\!\langle r_{n,1},\ldots,r_{n,n+u}\rangle\!\rangle,
$$
where the denominator is the closed normal subgroup generated by the relations. For a profinite group $G$ and a finite group $H$, we write $\Hom(G,H)$ for the set of continuous homomorphisms $G\to H$, and $\Sur(G,H)$ for its subset of surjections. For a random profinite group $X$ and a finite group $H$, we call $\bE(\#\Sur(X,H))$ the \emph{moment} of $X$ at $H$. Convergence in distribution of random profinite groups means weak convergence of their laws with respect to the topology defined in Section \ref{Sub22}.

\begin{thm}[{\cite[Theorem 1.1]{LW20}}]\label{thm1a}
Let $r_{n,1},\ldots,r_{n,n+u}$ be independent Haar-random elements of $\Fhat_n$. There exists a probability measure $\mu_u$ on isomorphism classes of profinite groups such that $\mc G_{n,u}$ converges in distribution to $\mu_u$ as $n\to\infty$.
\end{thm}

The measure $\mu_u$ is a nonabelian analogue of the distributions on abelian groups arising in the Cohen--Lenstra heuristics \cite{CL84}: the abelianization of a $\mu_u$-random group has the limiting distribution of the cokernels of Haar-random $n\times(n+u)$ matrices over $\widehat{\bZ}$ \cite[Example 12.4]{LW20}. Pagano and Shusterman \cite[Theorem 1.3]{PS26} proved that a profinite group with law $\mu_u$ is almost surely topologically finitely presented. Related random profinite group models, incorporating a Galois action and additional arithmetic restrictions, are used to predict distributions of Galois groups of maximal unramified extensions in \cite{LWZ24}.

Liu and Wood asked whether $\mu_u$ also arises from non-Haar distributions of the relations \cite[Section 1]{LW20}. They showed that $\mu_u$ arises from their model of random words of length at most $\ell$, in the iterated limit as $\ell$ and then $n$ tend to infinity \cite[Proposition 14.1]{LW20}. This raises the question of how the relation length must grow with $n$ for $\mc G_{n,u}$ to converge to $\mu_u$, and what happens at the threshold.

For abelian groups, the corresponding question is the universality problem for cokernels of random integral matrices. Wood \cite{Woo19} and Nguyen and Wood \cite{NW22} proved Cohen--Lenstra type limiting laws for matrices with independent entries satisfying suitable quantitative balance conditions. Their results distinguish convergence at a fixed finite set of primes from convergence of the full integral cokernel distribution. We study convergence to $\mu_u$ in three models of random relations: products of coset-balanced factors, random words of fixed length, and randomly ordered products of independent integer powers.

\subsection{Main results}\label{Sub11}

Fix $0<\epsilon\le1/2$. A probability measure $\nu$ on a profinite group $G$ is \emph{$\epsilon$-coset-balanced} if
$$
\nu(gK)\le1-\epsilon
\quad\text{for every proper open subgroup }K<G\text{ and every }g\in G.
$$
Under this hypothesis, Gorokhovsky proved a contraction estimate on each finite quotient of $G$, with a bound depending only on the quotient and $\epsilon$ \cite[Lemma 3.5]{Gor26}; see Proposition \ref{prop2c}. For positive integers $m_{n,j}$, $j\in[n+u]:=\{1,\ldots,n+u\}$, take mutually independent $\Fhat_n$-valued random variables $Y_{n,j,l}$ and set
$$
r_{n,j}:=Y_{n,j,1}\cdots Y_{n,j,m_{n,j}} \quad\text{for }j\in[n+u].
$$

\begin{thm}\label{thm1b}
\begin{enumerate}
\item Suppose that every $Y_{n,j,l}$ has an $\epsilon$-coset-balanced law. If
\begin{equation}\label{eq1a}
\frac{\min_{1\le j\le n+u}m_{n,j}}{\log n}\to\infty,
\end{equation}
then $\bE(\#\Sur(\mc G_{n,u},H))\to|H|^{-u}$ for every finite group $H$, and $\mc G_{n,u}\dist\mu_u$.

\item For every sequence of positive integers $q_n$ with $q_n/\log n\not\to\infty$, there exist $\epsilon$-coset-balanced probability measures $\nu_n$ on $\Fhat_n$ such that $\mc G_{n,u}$ does not converge weakly to $\mu_u$ when each relation has law $\nu_n^{*q_n}$.
\end{enumerate}
\end{thm}

For each fixed finite group $H$, Proposition \ref{prop2c} bounds the absolute relative error in the return probability for a product of $m$ factors by $|H|e^{-cm}$, where $c>0$ depends only on $H$ and $\epsilon$. Under \eqref{eq1a}, this bound is $o(n^{-1})$ uniformly over the relations and surjections, so the moment limits follow from \eqref{eq2n1}. Lemmas \ref{lem3a} and \ref{lem3b} give counterexamples for every fixed positive multiple of $\log n$, using coset-balanced laws on finite simple groups with arbitrarily small spectral gaps. The sharpness assertion is uniform over factor laws; particular laws may mix faster.

We next consider relations given by random words of fixed length in the generators and their inverses.
Let $\chi:\Fhat_n\to\bZ/2\bZ$ be the homomorphism sending every generator to $1$.
The uniform law on $\{x_1^{\pm1},\ldots,x_n^{\pm1}\}$ is supported on the nontrivial coset of $\ker\chi$, so it is not $\epsilon$-coset-balanced for any $\epsilon>0$. Every word of even length lies in $\ker\chi$, which explains the parity obstruction in Theorem \ref{thm1c}(d). For a positive integer $\ell_n$, set
$$
R_n:=x_{I_1}^{\varepsilon_1}\cdots x_{I_{\ell_n}}^{\varepsilon_{\ell_n}},
$$
where the indices $I_l$ are uniform on $[n]$, the signs $\varepsilon_l$ are uniform on $\{\pm1\}$, and all variables are independent. Let $G_{n,u}$ be the quotient of $\Fhat_n$ by $n+u$ independent copies of $R_n$. The length $\ell_n$ is counted before free reduction.

For odd lengths, the threshold is determined by the generators that occur in no relation. A fixed generator is absent from all $n+u$ relations with probability $(1-1/n)^{\ell_n(n+u)}$. Thus, when $\ell_n=\log n+O(1)$, the expected number of unused generators is $ne^{-\ell_n}(1+o(1))=e^{-(\ell_n-\log n)}(1+o(1))$. Each unused generator contributes a free factor to the quotient. For $s\in\bR$, set $\lambda_s:=e^{-s}$, take independent random variables $Z_s\sim\Pois(\lambda_s)$ and $X_w\sim\mu_w$, $w\in\bZ_{\ge0}$, and set
$$
\mc Y_{u,s}:=\Fhat_{Z_s}\amalg X_{u+Z_s}.
$$
Here $\amalg$ denotes the free profinite product \cite[Chapter 9]{RZ10}, and $\Fhat_0=1$.

\begin{thm}\label{thm1c}
Set $s_n:=\ell_n-\log n$.
\begin{enumerate}
\item $G_{n,u}\dist\mu_u$ if and only if $\ell_n$ is odd for all sufficiently large $n$ and $s_n\to\infty$.
\item If $n_j\to\infty$, $\ell_{n_j}$ is odd for every $j$, and $s_{n_j}\to s\in\bR$, then $G_{n_j,u}\dist\mc Y_{u,s}$.
\item If $s_n\to-\infty$, then the laws of $G_{n,u}$ have no weakly convergent subsequence.
\item If $\ell_n$ is even, then $G_{n,u}$ surjects onto $\bZ/2\bZ$ with probability $1$.
\end{enumerate}
\end{thm}

Along a subsequence as in Theorem \ref{thm1c}(b), deleting $r$ unused generators splits off a free factor $\Fhat_r$ and leaves $n+u$ relations on $n-r$ generators. For each fixed $r$, conditioned on every remaining generator occurring, the quotient on the remaining $n-r$ generators converges to $\mu_{u+r}$ (Proposition \ref{prop4g}).

The critical limit in Theorem \ref{thm1c}(b) is stated along subsequences because $\ell_n$ is integral. Indeed, if $s_n\to s\in\bR$ along all positive integers $n$, then
$$
\ell_{n+1}-\ell_n=\log(1+1/n)+s_{n+1}-s_n\to0,
$$
so $\ell_n$ would eventually be constant, contrary to $\ell_n\sim\log n$. Critical subsequences do occur: for example, $n_j=\lfloor e^{2j+1-s}\rfloor$ and $\ell_{n_j}=2j+1$ for all sufficiently large $j$ give $s_{n_j}\to s$.

We finally consider relations obtained by assigning an independent random integer exponent to each generator and multiplying the resulting powers in a uniformly random order.
For each $n$, take integer-valued random variables $\xi_{i,j}$, $i\in[n]$, $j\in[n+u]$, and uniform permutations $\sigma_{n,j}\in S_n$, all mutually independent. The variables $\xi_{i,j}$ may vary with $n$; we omit $n$ from the notation. Set
\begin{equation}\label{eq1c}
r^{\mr{ord}}_{n,j}:=x_{\sigma_{n,j}(1)}^{\xi_{\sigma_{n,j}(1),j}}\cdots x_{\sigma_{n,j}(n)}^{\xi_{\sigma_{n,j}(n),j}}, \quad G^{\mr{ord}}_{n,u}:=\Fhat_n/\langle\!\langle r^{\mr{ord}}_{n,1},\ldots,r^{\mr{ord}}_{n,n+u}\rangle\!\rangle.
\end{equation}
In a nonabelian quotient, the image of $r^{\mr{ord}}_{n,j}$ may depend on the order of the factors; the product is the endpoint of a random walk that uses each factor exactly once in a uniformly random order. For a prime $p$ and $0<\alpha\le1/2$, we say that an integer-valued random variable $\xi$ is \emph{$\alpha$-balanced modulo $p$} if
$$
\max_{a\in\bZ/p\bZ}\bP(\xi\equiv a\pmod p)\le1-\alpha,
$$
and that $\xi$ is \emph{$\alpha$-balanced} if this holds for every prime $p$, as in \cite{NW22}.

\begin{thm}\label{thm1d}
Suppose that, for each prime $p$, all $\xi_{i,j}$ are $\alpha_{n,p}$-balanced modulo $p$, with
\begin{equation}\label{eq1d}
\frac{n\alpha_{n,p}}{\log n}\to\infty.
\end{equation}
Then $\bE(\#\Sur(G^{\mr{ord}}_{n,u},H))\to|H|^{-u}$ for every finite group $H$, and $G^{\mr{ord}}_{n,u}\dist\mu_u$.
\end{thm}

Condition \eqref{eq1d} is imposed separately for each fixed prime, and the distributions of $\xi_{i,j}$ may vary with $n$, $i$, and $j$.

Set $\Xi_n:=(\xi_{i,j})$ and $B_n:=\bZ^n/\Xi_n\bZ^{n+u}$. For a profinite group $G$, write $G^{\mr{ab}}:=G/\overline{[G,G]}$ for its profinite abelianization, where $\overline{[G,G]}$ denotes the closure of the commutator subgroup of $G$. For an abelian group $B$, write $\widehat B$ for its profinite completion. 
Identify $\Fhat_n^{\mr{ab}}$ with $\widehat{\bZ}^n$ by sending $x_i$ to the standard basis vector $e_i$. Multiplication then becomes addition, so the image of $r^{\mr{ord}}_{n,j}$ is $\sum_{i=1}^n\xi_{i,j}e_i$, the $j$-th column of $\Xi_n$. Thus $(G^{\mr{ord}}_{n,u})^{\mr{ab}}\cong\widehat{\bZ}^n/\Xi_n\widehat{\bZ}^{n+u}\cong\widehat{B_n}$.

For a finite set $P$ of primes, write $B_n(P):=\prod_{p\in P}(B_n\otimes_{\bZ}\bZ_p)$. Theorem \ref{thm1d} yields
$$
\lim_{n\to\infty}\bP(B_n(P)\cong A)
=\frac1{|A|^u|\Aut(A)|}
\prod_{p\in P}\prod_{k=u+1}^{\infty}(1-p^{-k})
$$
for every finite abelian group $A$ whose order is divisible only by primes in $P$ (see Corollary \ref{cor6a}). This recovers the finite-prime cokernel law of Nguyen and Wood \cite[Theorem 4.1]{NW22}.

Gorokhovsky \cite[Theorem 3.1]{Gor23} studied random relations given by time-inhomogeneous random walks with symmetric step laws having positive mass at the identity and generating support, and obtained convergence to $\mu_u$ under the condition $\ell_n\epsilon_n/\log n\to\infty$, where $\ell_n$ is the number of steps and $\epsilon_n$ is the infimum of the positive step masses. For the lazy step law
$$
\nu_n:=\frac12\delta_1+\frac1{4n}\sum_{i=1}^n(\delta_{x_i}+\delta_{x_i^{-1}}),
\quad \epsilon_n=\frac1{4n},
$$
the condition becomes $\ell_n/(n\log n)\to\infty$. For non-lazy words of fixed odd length, Theorem \ref{thm1c} gives the sharp threshold $\ell_n-\log n\to\infty$ and the critical limits along subsequences.

The second author \cite[Theorem 1.3]{Lee26} proved universality of the cokernel distribution for square random matrices over the finite field $\bb F_p$ whose entries are independent and take any given value with probability at most $1-c\log n/n$, where $c>1$ is a fixed constant. For square matrices over the ring of $p$-adic integers $\bZ_p$, the authors and Yu \cite[Theorem 1.5]{JLY26} proved universality for $\alpha_n$-balanced entries with $\alpha_n=c\log n/n$ in the non-symmetric, symmetric and alternating models. Condition \eqref{eq1d} is stronger than what is needed in these matrix models, and we do not know whether $\liminf_{n\to\infty}n\alpha_n/\log n>1$ already suffices in the model \eqref{eq1c} for arbitrary $\alpha_n$-balanced exponents.

The paper is organized as follows. Section \ref{Sec2} collects the preliminaries, including Propositions \ref{prop2a} and \ref{prop2b}, which reduce convergence to $\mu_u$ to the moment limits $\bE(\#\Sur(\mc G_{n,u},H))\to|H|^{-u}$ for every finite group $H$. Independence of the relations expresses these moments in terms of return probabilities on $H$. Theorems \ref{thm1b} and \ref{thm1c} are proved in Sections \ref{Sec3} and \ref{Sec4}, respectively. For the critical limit in Theorem \ref{thm1c}(b), Proposition \ref{prop4g} establishes convergence conditional on every generator occurring in the unreduced sampled words. Section \ref{Sec5} proves Theorem \ref{thm1d}, and Section \ref{Sec6} recovers \cite[Theorem 4.1]{NW22} by abelianization.

\section{Preliminaries}\label{Sec2}

\subsection{Notation and conventions}\label{Sub21}
We use the following notation and conventions.
\begin{itemize}
\item Throughout the paper, $G$ denotes a profinite group, $H$ a finite group, and $u\in\bZ_{\ge0}$ is fixed. Asymptotic notation refers to $n\to\infty$ unless stated otherwise.

\item We use $\bP$ for probability, $\bE$ for expectation, and $\mc L(X)$ for the law of a random variable $X$. The symbols $\dist$ and $\xrightarrow{\bP}$ denote convergence in distribution and in probability, respectively. We write $\Pois(\lambda)$ for the Poisson distribution with mean $\lambda$, and $\Bin(n,p)$ for the binomial distribution with parameters $n$ and $p \in [0,1]$. For a nonempty finite set $A$, let $u_A$ be the uniform probability measure on $A$, $\delta_a$ the point mass at $a\in A$, and $\|\cdot\|_1$ the $\ell^1$-norm on $A$. For a subset or event $E$, let $\mb1_E$ denote its indicator.

\item For $B\subseteq G$, let $\langle B\rangle$ and $\langle\!\langle B\rangle\!\rangle$ denote the closed subgroup and closed normal subgroup generated by $B$, respectively. Bold lowercase letters denote tuples, such as $\mb h=(h_1,\ldots,h_n)\in H^n$ and $\mb a=(a_1,\ldots,a_N)\in A^N$.
For $\mb h \in H^n$, let $\varphi_{\mb h}:\Fhat_n\to H$ be the continuous homomorphism determined by $x_i\mapsto h_i$, and set
$$
\Gen_n(H):=\{\mb h\in H^n:\langle h_1,\ldots,h_n\rangle=H\}.
$$
The map $\mb h\mapsto\varphi_{\mb h}$ identifies $H^n$ with $\Hom(\Fhat_n,H)$ and $\Gen_n(H)$ with $\Sur(\Fhat_n,H)$. Every tuple in $H^n\setminus\Gen_n(H)$ lies in $K^n$ for some $K<H$. Thus, for every fixed finite group $H$,
$$
1-\frac{\#\Gen_n(H)}{|H|^n}
\le\sum_{K<H}\left(\frac{|K|}{|H|}\right)^n
\le\sum_{K<H}\frac{1}{2^n} \to0,
$$
which gives $\#\Gen_n(H) / |H|^n \to 1$ as $n \to \infty$.

\item We set $[n]:=\{1,\ldots,n\}$ for $n\ge0$ and $C_m:=\bZ/m\bZ$ for $m\ge1$. Unless stated otherwise, the letters $c$ and $C$ denote positive constants depending only on $H$ and the fixed parameters of the statement, and may change from one occurrence to the next.
\end{itemize}

For independent random relations $r_{n,1},\ldots,r_{n,n+u} \in \Fhat_n$, a surjection $\varphi:\Fhat_n\to H$ factors through $\mc G_{n,u}$ if and only if $\varphi(r_{n,j})=1$ for every $j$. Consequently,
\begin{equation}\label{eq2n1}
\begin{aligned}
\bE(\#\Sur(\mc G_{n,u},H))
&=\sum_{\varphi\in\Sur(\Fhat_n,H)}\bP(\varphi(r_{n,j})=1 \text{ for all }j\in[n+u])\\
&=\sum_{\varphi\in\Sur(\Fhat_n,H)}\prod_{j=1}^{n+u}\bP(\varphi(r_{n,j})=1)\\
&=\sum_{\mb h\in\Gen_n(H)}\prod_{j=1}^{n+u}\bP(\varphi_{\mb h}(r_{n,j})=1).
\end{aligned}
\end{equation}
The second equality uses independence of the relations, and the
last equality uses the identification of $\Sur(\Fhat_n,H)$ with $\Gen_n(H)$.

\subsection{Moments and convergence in distribution}\label{Sub22}

For a finite set $\mc S$ of isomorphism classes of finite groups, let $\mc S^{\mr{cl}}$ be its closure under taking subgroups, quotients, and finite direct products. Set
$$
G^{\mc S}:=G/\bigcap_{\substack{A\in\mc S^{\mr{cl}}\\ \varphi\in\Hom(G,A)}}\ker\varphi,
$$
the pro-$\mc S^{\mr{cl}}$ completion of $G$. Let $\mc P$ be the set of isomorphism classes of profinite groups $G$ such that $G^{\mc S}$ is finite for every finite $\mc S$. We equip $\mc P$ with the topology generated by the clopen sets \cite[Sections 2.4 and 3]{LW20}
$$
U_{\mc S,A}:=\{G\in\mc P:G^{\mc S}\cong A\},
$$
where $\mc S$ ranges over all finite sets of isomorphism classes of finite groups and $A$ ranges over $\mc S^{\mr{cl}}$. Every finitely generated profinite group belongs to $\mc P$ \cite[Corollary 15.72]{Neu67}. The counts $\#\Hom(G,H)$ and $\#\Sur(G,H)$ depend only on $G^{\{H\}}$ and are finite and locally constant on $\mc P$.

\begin{prop}[{\cite[Corollary 2.22]{Woo22}}]\label{prop2a}
Let $X_n,Y$ be random groups in $\mc P$. Suppose that there exist constants $c,C>0$, independent of $H$, such that
$$
\lim_{n\to\infty}\bE(\#\Sur(X_n,H))
=\bE(\#\Sur(Y,H))\le C|H|^c
\quad\text{for every finite group }H.
$$
Then $X_n\dist Y$.
\end{prop}

\begin{prop}[{\cite[Theorem 1.1]{LW20} and \cite[Lemma 4.1]{Lee24}}]\label{prop2b}
Fix $w\in\bZ_{\ge0}$ and let $X_w\sim\mu_w$. For $n\ge1$, let $Q_{n,w}$ be the quotient of $\Fhat_n$ by $n+w$ independent Haar-random relations. Then $Q_{n,w}\dist\mu_w$ and $\bE(\#\Sur(X_w,H))=|H|^{-w}$ for every finite $H$.
\end{prop}

\subsection{Convolution estimates}\label{Sub23}

We equip $\ell^2(H)$ with counting measure and write $\ell^2_0(H)$ for its subspace of functions with zero sum. We regard a probability measure $\nu$ on a subgroup $K\le H$ as a probability measure on $H$ by setting $\nu(h)=0$ for every $h\in H\setminus K$. For probability measures $\nu,\eta$ on $H$, convolution and inversion are defined by
$$
(\nu*\eta)(x):=\sum_{g\in H}\nu(g)\eta(g^{-1}x),
\quad \check\nu(x):=\nu(x^{-1}).
$$
For $m\in\bZ_{\ge0}$, $\nu^{*m}$ denotes the $m$-fold convolution power, with $\nu^{*0}:=\delta_1$. We write $\bar\nu:=(\nu+\check\nu)/2$ for the symmetrization of $\nu$, and we call $\nu$ symmetric if $\check\nu=\nu$.

Set $R_gf(x):=f(xg)$ and $T_\nu:=\sum_{g\in H}\nu(g)R_g$. These operators fix constants, preserve $\ell^2_0(H)$, and satisfy
$$
T_{\nu*\eta}=T_\nu T_\eta,\quad
T_{\check\nu}=T_\nu^*,\quad
\Tr(T_\nu)=|H|\nu(1).
$$
Right convolution is given by $f*\nu=T_{\check\nu}f$. For a linear operator $T$ on a subspace $V$ of $\ell^2(H)$, write $\|T\|_\op:=\sup\{\|Tf\|_2:f\in V,\ \|f\|_2\le1\}$ for its operator norm.

\begin{prop}[{\cite[Lemma 3.5]{Gor26}}]\label{prop2c}
Let $0<\epsilon\le1/2$. For every $\epsilon$-coset-balanced probability measure $\nu$ on $H$,
\begin{equation}\label{eq2a}
\|T_\nu|_{\ell^2_0(H)}\|_\op
\le\exp\left(-\frac{\epsilon}{2|H|^3}\right).
\end{equation}
\end{prop}

\subsection{Types and entropy estimates}\label{Sub24}

For a nonempty finite set $A$, an integer $N\ge1$, and $\mb a\in A^N$, let
$$
v_{\mb a}:=\frac1N\sum_{i=1}^N\delta_{a_i}
$$
be the \emph{type} of $\mb a$. For a probability vector $v$ on $A$, write $\supp(v):=\{a\in A:v(a)>0\}$ for its support and $\mr H(v):=-\sum_{a\in A}v(a)\log v(a)$ for its \emph{entropy}, with the conventions $0\log0=0$ and $\log0=-\infty$. There are at most $(N+1)^{|A|}$ types of tuples in $A^N$, and each type $v$ is realized by at most $e^{N\mr H(v)}$ tuples. Pinsker's inequality \cite[Chapter 11]{CT06} gives 
$$
\mr H(v)\le\log|A|-\frac12\|v-u_A\|_1^2.
$$
Hence, for every $\theta\ge0$, summing over the types satisfying $\|v-u_A\|_1>\theta$ gives
\begin{equation}\label{eq2b}
\#\{\mb a\in A^N:\|v_{\mb a}-u_A\|_1>\theta\}
\le (N+1)^{|A|}|A|^N e^{-N\theta^2/2}.
\end{equation}
For probability vectors on a fixed finite set, $v_n\to v$ means $\|v_n-v\|_1\to0$.

\begin{lem}\label{lem2d}
For each $n\ge1$, let $Q_n:\{v_{\mb h}:\mb h\in H^n\}\to[0,1]$ be a function. Assume that
$$
\limsup_{n\to\infty}Q_n(v_n)
\le|\langle\supp(v)\rangle|^{-1}
$$
whenever $v_n\to v$, where $v_n=v_{\mb h_n}$ with $\mb h_n\in H^n$ and $v$ is a probability vector on $H$.
For each $K\le H$, let $\mc N_K$ be an $\ell^1$-neighborhood of $u_K$ in the space of probability vectors on $H$. Then there exists $c>0$ such that, for all sufficiently large $n$,
$$
\mr H(v)+\log Q_n(v)\le-c
$$
for every type $v$ of a tuple in $H^n$ lying outside $\bigcup_{K\le H}\mc N_K$.
\end{lem}

\begin{proof}
Suppose that the conclusion fails. Then there exist $n_j\to\infty$ and types $v_j$ of tuples in $H^{n_j}$, each outside the prescribed neighborhoods, such that $\liminf_j(\mr H(v_j)+\log Q_{n_j}(v_j))\ge0$. Since the space of probability vectors on $H$ is compact in the $\ell^1$-topology, we may pass to a subsequence and assume that $v_j\to v$ for some probability vector $v$. 

Set $K:=\langle\supp(v)\rangle$. The hypothesis, continuity of entropy, and Pinsker's inequality on $K$ give
$$
0\le\limsup_j(\mr H(v_j)+\log Q_{n_j}(v_j))
\le\mr H(v)-\log|K|
\le-\frac12\|v-u_K\|_1^2.
$$
Thus $v=u_K$, so $v_j\in\mc N_K$ for all sufficiently large $j$, contrary to the choice of $v_j$.
\end{proof}

\section{Coset-balanced products}\label{Sec3}

We first construct coset-balanced measures with arbitrarily slow mixing, which will be used in the proof of Theorem \ref{thm1b}(b).

\begin{lem}\label{lem3a}
For every $\beta_0\in(0,1)$, there exist a finite nonabelian simple group $H$ and a symmetric $1/2$-coset-balanced probability measure $\zeta$ on $H$ such that
\begin{equation}\label{eq3a}
|H|\zeta^{*m}(1)\ge1+\beta_0^m
\quad\text{for every }m\in\bZ_{\ge1}.
\end{equation}
\end{lem}

\begin{proof}
Let $d\ge5$ be odd, and let $H=A_d$ be the alternating group on $d$ letters acting on $\bZ/d\bZ$. Consider a $d$-cycle $a:=(0\ 1\ \cdots\ d-1)$ and a $3$-cycle $b:=(-1\ 0\ 1)$ in $H$. The conjugates $a^iba^{-i}=(i-1\ i\ i+1)$ ($i\in\bZ/d\bZ$) are the consecutive $3$-cycles, which generate $H$. Hence $\langle a,b\rangle=H$. 
Let $\eta$ be the uniform measure on $\{1,a,b,ba\}$. If a coset $gK$ of a subgroup $K\le H$ contains three of these four elements, then $a,b\in K$, hence $K=H$. Hence $\eta$ is $1/2$-coset-balanced. The symmetric measure $\zeta:=\check\eta*\eta$ is again $1/2$-coset-balanced, because 
$$
\zeta(gK)=\sum_{x\in H}\eta(x)\eta(xgK)\le\frac12\sum_{x\in H}\eta(x)=\frac12
$$ 
for every proper subgroup $K<H$ and every $g\in H$. Moreover, $T_\zeta=T_\eta^*T_\eta$ is positive semidefinite.

Set $A:=\{2,\ldots,(d+1)/2\}\subseteq\bZ/d\bZ$, $\rho:=|A|/d=(d-1)/(2d)$, and $f(g):=\mb1_A(g^{-1}(0))-\rho$. Since $g^{-1}(0)$ is uniformly distributed on $\bZ/d\bZ$ as $g$ ranges over $H$, we have $f\in\ell^2_0(H)$ and $\|f\|_2^2=|H|\rho(1-\rho)$. 
The permutation $b$ fixes every point of $A$, so $R_bf=R_{b^{-1}}f=f$. Also, $A$ and its translate by $a$ differ at exactly two points, so $\|R_af-f\|_2^2=2|H|/d$. Since $R_a$ is unitary and $f$ is real-valued,
$$
\langle R_af,f\rangle
=\|f\|_2^2-\frac12\|R_af-f\|_2^2
=\|f\|_2^2-\frac{|H|}{d}.
$$
Since
$$
T_\eta=\frac14(R_1+R_a+R_b+R_{ba})=\frac14(I+R_b)(I+R_a)
$$ 
and $R_b^*f=R_{b^{-1}}f=f$, we have
$$
\langle T_\eta f,f\rangle
=\frac14\langle (I+R_a)f,(I+R_b^*)f\rangle
=\frac12(\|f\|_2^2+\langle R_af,f\rangle)
=\|f\|_2^2-\frac{|H|}{2d}.
$$

Let $\beta_d$ be the largest eigenvalue of $T_\zeta|_{\ell^2_0(H)}$. Since $T_\zeta=T_\eta^*T_\eta$, we have $\beta_d \ge \|T_\eta f\|_2^2 / \|f\|_2^2$. Moreover, the inequality $\langle T_\eta f,f\rangle^2 \le \|T_\eta f\|_2^2 \|f\|_2^2$ gives
$$
\frac{\|T_\eta f\|_2^2}{\|f\|_2^2}
\ge \frac{\langle T_\eta f,f\rangle^2}{\|f\|_2^4}
= \left(1-\frac{|H|}{2d\|f\|_2^2}\right)^2
= \left(1-\frac1{2d\rho(1-\rho)}\right)^2.
$$
By \eqref{eq2a}, $\beta_d\le\exp(-1/(4|H|^3))<1$. Hence $\beta_d\to1$ as $d\to\infty$, and we may choose $d$ so that $\beta_d>\beta_0$. Since $T_\zeta$ is positive semidefinite,
\begin{equation*}
|H|\zeta^{*m}(1)=\Tr(T_\zeta^m)\ge1+\beta_0^m. \qedhere
\end{equation*}
\end{proof}

The next lemma lifts the measure from Lemma \ref{lem3a} to $\Fhat_n$ while preserving coset balance.

\begin{lem}\label{lem3b}
Let $H$ be a finite nonabelian simple group, and let $\zeta$ be a $1/2$-coset-balanced probability measure on $H$. For all sufficiently large $n$, choose representatives $\varphi_1,\ldots,\varphi_{M_n}$ of the $\Aut(H)$-orbits in $\Sur(\Fhat_n,H)$, so that $M_n=\#\Gen_n(H)/|\Aut(H)|$. Then there exists a $1/2$-coset-balanced probability measure $\nu_n$ on $\Fhat_n$ such that
$$
(\varphi_1,\ldots,\varphi_{M_n})_*\nu_n=\zeta^{\otimes M_n}.
$$
\end{lem}

\begin{proof}
Consider $\Phi:=(\varphi_1,\ldots,\varphi_{M_n}):\Fhat_n\to H^{M_n}$. Its image is a subdirect subgroup of $H^{M_n}$, since each $\varphi_i$ is surjective. If $\Phi$ were not surjective, then \cite[Theorem 4.16(iii)]{PS18} implies that $\varphi_j=\alpha\circ\varphi_i$ for some $i\ne j$ and some $\alpha\in\Aut(H)$, which contradicts the choice of representatives. Thus $\Phi$ is surjective. The product measure $\zeta^{\otimes M_n}$ is $1/2$-coset-balanced by \cite[Lemma 3.4(2)]{Gor26}.

Lift it to a measure $\nu_n$ on $\Fhat_n$ by equipping each fiber of $\Phi$ with a translate of the normalized Haar measure on $\ker\Phi$. Let $K<\Fhat_n$ be a proper open subgroup and $g\in\Fhat_n$. If $\Phi(K)<H^{M_n}$, then
$$
\nu_n(gK)\le\zeta^{\otimes M_n}(\Phi(g)\Phi(K))\le\frac12.
$$
If $\Phi(K)=H^{M_n}$, each fiber meets $gK$ in a coset of $\ker\Phi\cap K$, so
\begin{equation*}
\nu_n(gK)=[\ker\Phi:\ker\Phi\cap K]^{-1}
=[\Fhat_n:K]^{-1}\le\frac12. \qedhere
\end{equation*}
\end{proof}

\begin{proof}[Proof of Theorem \ref{thm1b}]
\begin{enumerate}
\item Assume \eqref{eq1a}, and fix a finite group $H$. For $\varphi\in\Sur(\Fhat_n,H)$ and $j\in[n+u]$, set $\eta_l:=\varphi_*\mc L(Y_{n,j,l})$ for $1\le l\le m_{n,j}$. Each $\eta_l$ is $\epsilon$-coset-balanced, and $\mu:=\eta_1*\cdots*\eta_{m_{n,j}}$ is the law of $\varphi(r_{n,j})$. For $f:=\delta_1-u_H$, we obtain
\begin{align*}
|H|\bP(\varphi(r_{n,j})=1)-1
=|H|\left(\mu(1)-\frac1{|H|}\right)
=|H|\sum_{g\in H}\mu(g)f(g)
=|H|(T_{\eta_1}\cdots T_{\eta_{m_{n,j}}}f)(1).
\end{align*}
Since $f\in\ell^2_0(H)$ and $\|f\|_2^2=1-|H|^{-1}\le1$, Proposition \ref{prop2c} gives
\begin{align*}
| |H|\bP(\varphi(r_{n,j})=1)-1 |
&\le |H|\|T_{\eta_1}\cdots T_{\eta_{m_{n,j}}}f\|_2\\
&\le |H|\prod_{l=1}^{m_{n,j}}\|T_{\eta_l}|_{\ell^2_0(H)}\|_\op\\
&\le |H|\exp\left(-\frac{\epsilon m_{n,j}}{2|H|^3}\right).
\end{align*}
Taking the maximum over $j$, $\varphi$ and using
\eqref{eq1a}, we conclude that
$$
\rho_n:=\max_{j,\varphi}
\left||H|\bP(\varphi(r_{n,j})=1)-1\right|
\le |H|\exp\left(-\frac{\epsilon\min_jm_{n,j}}{2|H|^3}\right)
=o(n^{-1}).
$$

The relations are independent, so uniformly in $\varphi$,
$$
\prod_{j=1}^{n+u}\bP(\varphi(r_{n,j})=1)
=|H|^{-n-u}\exp(O(n\rho_n))
=|H|^{-n-u}(1+o(1)).
$$
By \eqref{eq2n1} and the preceding uniform estimate, we have
$$
\bE(\#\Sur(\mc G_{n,u},H))
=\frac{\#\Gen_n(H)}{|H|^{n+u}}(1+o(1)) \to |H|^{-u}.
$$
Propositions \ref{prop2a} and \ref{prop2b} give $\mc G_{n,u}\dist\mu_u$.

\item Suppose that $q_n/\log n\not\to\infty$. Then there exist $C>0$ and an infinite set $\mc N\subseteq\bZ_{\ge1}$ such that $q_n\le C\log n$ for all $n\in\mc N$. Fix $\beta_0\in(e^{-1/C},1)$, and let $H$ and $\zeta$ be as in Lemma \ref{lem3a}. For all sufficiently large $n$, take $\nu_n$ as in Lemma \ref{lem3b}, and use Haar measure for the remaining $n$. Each $\nu_n$ is $\epsilon$-coset-balanced because $\epsilon\le1/2$. Give each relation the law $\nu_n^{*q_n}$.

Let $V_n:=\#\Sur(\mc G_{n,u},H)/|\Aut(H)|$. Since convolution commutes with pushforward, Lemma \ref{lem3b} gives
$$
(\varphi_1,\ldots,\varphi_{M_n})_*\nu_n^{*q_n}
=(\zeta^{*q_n})^{\otimes M_n}.
$$
The preceding identity and independence of the relations imply that the events $E_i:=\{\varphi_i(r_{n,j})=1\text{ for all }j\in[n+u]\}$ ($1 \le i \le M_n$) are independent, with
$$
\bP(E_i)
=\prod_{j=1}^{n+u}\bP(\varphi_i(r_{n,j})=1)
=\zeta^{*q_n}(1)^{n+u}.
$$
Hence $V_n=\sum_{i=1}^{M_n}\mb1_{E_i}\sim\Bin(M_n,\zeta^{*q_n}(1)^{n+u})$. By \eqref{eq3a} and the inequality $\log(1+x) \ge x/2$ for $x \in [0,1]$,
\begin{align*}
\bE(V_n)
&=\frac{M_n}{|H|^{n+u}}(|H|\zeta^{*q_n}(1))^{n+u} \\
&\ge\frac{1+o(1)}{|\Aut(H)||H|^u} (1+\beta_0^{q_n})^{n+u} \\
&\ge\frac{1+o(1)}{|\Aut(H)||H|^u} \exp\left(\frac{n\beta_0^{q_n}}2\right).
\end{align*}
For $n\in\mc N$, $n\beta_0^{q_n}\ge n^{1+C\log\beta_0}\to\infty$, since $\beta_0>e^{-1/C}$. Hence $\bE(V_n)\to\infty$ along $\mc N$. Since $V_n$ has a binomial distribution, $\Var(V_n)\le\bE(V_n)$. Chebyshev's inequality therefore implies that $V_n\xrightarrow{\bP}\infty$ along $\mc N$. On the other hand, if $\mc G_{n,u}$ converged weakly to $\mu_u$, then $(V_n)$ would be tight because $\#\Sur(\,\cdot\,,H)$ is locally constant on $\mc P$, a contradiction. \qedhere
\end{enumerate}
\end{proof}

\section{Fixed-length words}\label{Sec4}

For $w\in\bZ_{\ge0}$ and $\ell\in\bZ_{\ge1}$, let $G_{n,w}[\ell]$ be the quotient of $\Fhat_n$ by $n+w$ independent words of length $\ell$, sampled as in Section \ref{Sec1}. The random group $G_{n,u}$ defined in the introduction corresponds to the choice $w=u$ and $\ell=\ell_n$.

\subsection{Return probability estimates}\label{Sub41}

For $N\ge1$ and $\mb h\in H^N$, let $I$ and $\varepsilon$ be independent and uniform on $[N]$ and $\{\pm1\}$, respectively. The random letter $x_I^\varepsilon$ maps under $\varphi_{\mb h}$ to $h_I^\varepsilon$, whose law is
$$
\bar v_{\mb h}=\frac12(v_{\mb h}+\check v_{\mb h})=\frac1{2N}\sum_{i=1}^N(\delta_{h_i}+\delta_{h_i^{-1}}).
$$
Set $q_{\ell,H}(\mb h):=\bar v_{\mb h}^{*\ell}(1)$, the probability that the image of a random word of length $\ell$ under $\varphi_{\mb h}$ is the identity. Then \eqref{eq2n1} gives
\begin{equation}\label{eq4a}
\bE(\#\Sur(G_{n,w}[\ell],H)) =\sum_{\mb h\in\Gen_n(H)}q_{\ell,H}(\mb h)^{n+w}.
\end{equation}

\begin{lem}\label{lem4a}
Let $N\ge1$, $\mb h\in H^N$ and $\theta:=\|v_{\mb h}-u_H\|_1<1$. Then
\begin{equation}\label{eq4b}
\left|q_{\ell,H}(\mb h)-|H|^{-1}\right|\le\theta^\ell
\quad\text{for every }\ell\ge1.
\end{equation}
Moreover, if $\nu$ is a symmetric probability measure on $H$, $\ell$ is odd and $m$ is even with $2\le m\le\ell$, then
\begin{equation}\label{eq4c}
\nu^{*\ell}(1)\le\nu^{*m}(1).
\end{equation}
\end{lem}

\begin{proof}
Since $T_{u_H}$ vanishes on $\ell^2_0(H)$ and each $R_g$ is unitary, 
\begin{align*}
\|T_{\bar v_{\mb h}}f\|_2
&=\|(T_{\bar v_{\mb h}}-T_{u_H})f\|_2\\
&=\left\|\sum_{g\in H}(\bar v_{\mb h}(g)-u_H(g))R_gf\right\|_2\\
&\le\sum_{g\in H}|\bar v_{\mb h}(g)-u_H(g)|\|R_gf\|_2\\
&=\|\bar v_{\mb h}-u_H\|_1\|f\|_2
\end{align*}
for every $f \in \ell^2_0(H)$. Thus every eigenvalue of $T_{\bar v_{\mb h}}|_{\ell^2_0(H)}$ has absolute value at most $\|\bar v_{\mb h}-u_H\|_1$. Moreover, symmetrization gives $\|\bar v_{\mb h}-u_H\|_1\le\theta$. Writing the eigenvalues of $T_{\bar v_{\mb h}}$ as $1,\lambda_2,\ldots,\lambda_{|H|}$, the trace identity gives
$$
\left|q_{\ell,H}(\mb h)-|H|^{-1}\right|
=\frac1{|H|}\left|\sum_{j=2}^{|H|}\lambda_j^\ell\right|
\le\theta^\ell.
$$
For $\lambda\in[-1,1]$, we have $\lambda^\ell\le\lambda^m$ when $\ell$ is odd and $m\le\ell$ is even. Thus
\begin{equation*}
\nu^{*\ell}(1)=\frac1{|H|}\Tr(T_\nu^\ell)
\le\frac1{|H|}\Tr(T_\nu^m)=\nu^{*m}(1). \qedhere
\end{equation*}
\end{proof}

\begin{lem}\label{lem4b}
Let $(\ell_n)_{n\ge1}$ be a sequence of odd positive integers with $\ell_n\to\infty$. If $v_{\mb h_n}\to v$ for $\mb h_n\in H^n$, then
$$
\limsup_{n\to\infty}q_{\ell_n,H}(\mb h_n)
\le|\langle\supp(v)\rangle|^{-1}.
$$
\end{lem}

\begin{proof}
Set $L:=\langle\supp(v)\rangle$. For each $g \in H$, $R_g$ is unitary, so $\|R_gf\|_2=\|f\|_2$. If $T_{\bar v}f=f$, then
$$
\sum_{g\in H}\bar v(g)\|R_gf-f\|_2^2=2\|f\|_2^2-2\mr{Re}\langle T_{\bar v}f,f\rangle=0.
$$
Thus $R_gf=f$ for every $g\in\supp(\bar v)$, and hence for every $g\in L$. Conversely, every right $L$-invariant function is fixed by $T_{\bar v}$. Thus the $1$-eigenspace of $T_{\bar v}$ consists of the right $L$-invariant functions and has dimension $[H:L]$. 
By eigenvalue continuity, all but $[H:L]$ eigenvalues of $T_{\bar v_{\mb h_n}}$ are at most $1-c$ for some fixed $c\in(0,1)$ and all sufficiently large $n$. Hence
\begin{equation*}
q_{\ell_n,H}(\mb h_n)
=\frac1{|H|}\Tr(T_{\bar v_{\mb h_n}}^{\ell_n})
\le\frac{[H:L]}{|H|}+\frac{|H|-[H:L]}{|H|}(1-c)^{\ell_n}
\le|L|^{-1}+(1-c)^{\ell_n}. \qedhere
\end{equation*}
\end{proof}

For each fixed $K\le H$ and $\mb h\in H^n$, set $k(\mb h):=\#\{i:h_i\notin K\}$ and let $\mb h_K$ be the subtuple of $\mb h$ consisting of the coordinates that belong to $K$.

\begin{lem}\label{lem4c}
Let $K\le H$, $\mb h\in H^n$, $k:=k(\mb h)\le n/(8|H|)$, and $\|v_{\mb h_K}-u_K\|_1\le 1/(8|H|)$. For odd $\ell\ge1$ and $w\in\bZ_{\ge0}$, the following hold.
\begin{enumerate}
\item If $\ell k/n\le 1/(4|H|)$, then
\begin{equation}\label{eq4d}
q_{\ell,H}(\mb h)^{n+w}
\le|K|^{-n}\exp\left(-\ell k+|H|\frac{\ell^2k^2}{n}+\frac{2|H|n}{(8|H|)^\ell}\right).
\end{equation}
\item If $\ell k/n\ge 1/(4|H|)$, then
$$
q_{\ell,H}(\mb h)\le|K|^{-1}\left(1-\frac1{16|H|}\right).
$$
\end{enumerate}
\end{lem}

\begin{proof}
Write $\beta:=k/n$. In a word of length $t$, the number $B_t$ of letters whose images lie outside $K$ has distribution $\Bin(t,\beta)$. Conditional on $B_t=0$, the images are independent with law $\bar v_{\mb h_K}$. Thus the contribution of $B_t=0$ to $q_{t,H}(\mb h)$ is
$$
\bP(B_t=0)q_{t,K}(\mb h_K)=(1-\beta)^t q_{t,K}(\mb h_K).
$$
If $B_t=1$, the image cannot be $1$, since $K(H\setminus K)K\subseteq H\setminus K$. Hence \eqref{eq4b} and $\bP(B_t\ge2)\le\binom t2\beta^2$ give, for every $t\ge1$,
$$
q_{t,H}(\mb h)
\le(1-\beta)^t q_{t,K}(\mb h_K)+\binom t2\beta^2
\le|K|^{-1}e^{-t\beta}+(8|H|)^{-t}+\frac{(t\beta)^2}{2}.
$$
\begin{enumerate}
\item Put $x:=\ell\beta\le 1/(4|H|)$. Since $e^x\le2$,
$$
|K|q_{\ell,H}(\mb h)\le e^{-x}\left(1+|H|x^2+\frac{2|H|}{(8|H|)^\ell}\right).
$$
Using $q_{\ell,H}(\mb h)\le1$, $w\ge0$, and $(1+y)^n\le e^{ny}$, we obtain
\begin{align*}
q_{\ell,H}(\mb h)^{n+w}
&\le q_{\ell,H}(\mb h)^n\\
&\le|K|^{-n}e^{-nx}\left(1+|H|x^2+\frac{2|H|}{(8|H|)^\ell}\right)^n\\
&\le|K|^{-n}\exp\left(-\ell k+|H|\frac{\ell^2k^2}{n}+\frac{2|H|n}{(8|H|)^\ell}\right).
\end{align*}

\item Set $m:=2\lfloor 1/(8|H|\beta)\rfloor$. Since $\beta\le 1/(8|H|)$ and $\ell \beta \ge 1/(4|H|)$, we have $2\le m\le\ell$ and $1/(8|H|)\le m\beta\le 1/(4|H|)$. By \eqref{eq4c} and the preceding estimate with $t=m$, we have
\begin{align*}
q_{\ell,H}(\mb h)
&\le q_{m,H}(\mb h)\\
&\le|K|^{-1}e^{-m\beta}+(8|H|)^{-m}+\frac{(m\beta)^2}{2}\\
&\le|K|^{-1}\left(1-m\beta+\frac{(m\beta)^2}{2}\right)+m\beta \frac1{8|H|}+\frac{m\beta}{2}\frac1{4|H|}\\
&\le|K|^{-1}\left(1-\frac{m\beta}{2}\right)\\
&\le|K|^{-1}\left(1-\frac1{16|H|}\right). \qedhere
\end{align*}
\end{enumerate}
\end{proof}

\subsection{Moment estimates}\label{Sub42}

Fix $w\in\bZ_{\ge0}$ and an infinite $\mc N\subseteq\bZ_{\ge1}$ along which $\ell_n$ is odd and $\ell_n\ge\log n+b$ for some $b\in\bR$. Throughout this subsection, $n$ tends to infinity through $\mc N$. Let $\mc D_n$ be the event that every generator occurs in the unreduced relations defining $G_{n,w}[\ell_n]$. For finite $H$ and $\mb h\in H^n$, let $\mc A_{\mb h}$ be the event that all relations map to $1$ under $\varphi_{\mb h}$; then $\bP(\mc A_{\mb h})=q_{\ell_n,H}(\mb h)^{n+w}$.

\begin{prop}\label{prop4d}
Fix $w\in\bZ_{\ge0}$, and let $\mc N\subseteq\bZ_{\ge1}$ be infinite. Suppose that $\ell_n$ is odd for every $n\in\mc N$ and $\ell_n-\log n\to\infty$ as $n\to\infty$ through $\mc N$. Then, along $\mc N$,
$$
\bE(\#\Sur(G_{n,w}[\ell_n],H))\to|H|^{-w}
$$
for every finite group $H$.
\end{prop}

Fix a nontrivial finite group $H$. For a probability vector $v$ on $H$, set
$$
Q_n(v):=\bar v^{*\ell_n}(1),\quad
\Phi_n(v):=\mr H(v)+\log Q_n(v),\quad
\mc B_n:=\{\mb h\in H^n:\|v_{\mb h}-u_H\|_1\le1/(8|H|)\}.
$$
Then $Q_n(v_{\mb h})=q_{\ell_n,H}(\mb h)$, and every $\mb h\in\mc B_n$ generates $H$, since $v_{\mb h}(g)\ge7/(8|H|)>0$ for all $g\in H$. Now \eqref{eq2b} gives $\#\mc B_n=|H|^n(1+o(1))$, and \eqref{eq4b} gives a uniform estimate
$$
q_{\ell_n,H}(\mb h)^{n+w}=|H|^{-n-w}(1+o(1)) \quad (\mb h \in \mc B_n).
$$
In particular, we have
\begin{equation}\label{eq4f}
\lim_{n \to \infty} \sum_{\mb h\in\mc B_n}q_{\ell_n,H}(\mb h)^{n+w}=|H|^{-w}.
\end{equation}

\begin{prop}\label{prop4e}
\begin{enumerate}
\item For every $t>0$,
$$
\sum_{\substack{\mb h\in H^n\\ \Phi_n(v_{\mb h})<-t}}q_{\ell_n,H}(\mb h)^{n+w}
\le(n+1)^{|H|}e^{-nt}.
$$
\item For all sufficiently large $n$, every $\mb h\in\Gen_n(H)\setminus\mc B_n$ with $\Phi_n(v_{\mb h})\ge-n^{-2/3}$ admits $K<H$ such that
\begin{equation} \label{eq4n1}
1\le k(\mb h)\le n/(|H|\ell_n^2),\quad
\|v_{\mb h_K}-u_K\|_1\le1/(8|H|).
\end{equation}
\end{enumerate}
\end{prop}

\begin{proof}
\begin{enumerate}
\item Since there are at most $(n+1)^{|H|}$ types and each type $v$ occurs at most $e^{n\mr H(v)}$ times,
$$
\sum_{\Phi_n(v_{\mb h})<-t}q_{\ell_n,H}(\mb h)^{n+w}
\le\sum_{\substack{v\text{ a type}\\ \Phi_n(v)<-t}}e^{n\mr H(v)}Q_n(v)^n
\le(n+1)^{|H|}e^{-nt}.
$$

\item Suppose that (b) fails along a subsequence, and choose a counterexample $\mb h_n$ for each $n$ in this subsequence. Set $v_n:=v_{\mb h_n}$. Since the space of probability vectors on $H$ is compact, after passing to a further subsequence we may assume that $v_n\to v$ for some probability vector $v$ on $H$. All limits below are along this subsequence.

Set $K:=\langle\supp(v)\rangle$. Since $\Phi_n(v_n)\ge-n^{-2/3}$, Lemma \ref{lem4b}, continuity of entropy, and Pinsker's inequality on $K$ give
$$
0\le\limsup_{n\to\infty}\Phi_n(v_n)
\le\mr H(v)-\log|K|
\le-\frac12\|v-u_K\|_1^2.
$$
Thus $v=u_K$. Moreover, $\mb h_n\notin\mc B_n$ implies $\|v-u_H\|_1\ge1/(8|H|)$, so $K<H$.

Set $k_n:=k(\mb h_n)$ and $v_n(K):=\sum_{a\in K}v_n(a)$. Then $v_n\to u_K$ gives 
$$
\frac{k_n}{n}= 1-v_n(K) \to 0, \quad
\|v_{(\mb h_n)_K}-u_K\|_1=\sum_{a\in K}\left|\frac{v_n(a)}{v_n(K)}-\frac1{|K|}\right|\to0,
$$
so the hypotheses of Lemma \ref{lem4c} hold for all sufficiently large $n$. Since $\mb h_n$ generates $H$, we have $k_n\ge1$. The condition \eqref{eq4n1} fails for $\mb h = \mb h_n$, hence $k_n > n/(|H|\ell_n^2)$.
If $\ell_n k_n/n\ge1/(4|H|)$ for infinitely many $n$, then Lemma \ref{lem4c}(b) and the continuity of entropy give, along a further subsequence,
$$
-n^{-2/3} \le\Phi_n(v_n) \le\mr H(v_n)-\log|K|+\log\left(1-\frac1{16|H|}\right)
=\log\left(1-\frac1{16|H|}\right)+o(1),
$$
which is a contradiction. Hence $\ell_n k_n/n<1/(4|H|)$ for all sufficiently large $n$.

Splitting the entropy between $K$ and its complement gives
\begin{align*}
\mr H(v_n)
&\le\left(1-\frac{k_n}{n}\right)\log\frac{|K|}{1-k_n/n}+\frac{k_n}{n}\log\frac{n(|H|-|K|)}{k_n}\\
&=\log|K|+\frac{k_n}{n}\log\frac{n(|H|-|K|)}{|K|k_n}-\left(1-\frac{k_n}{n}\right)\log\left(1-\frac{k_n}{n}\right)\\
&\le\log|K|+\frac{k_n}{n}\log\frac{n(|H|-|K|)}{|K|k_n}+\frac{k_n}{n}\\
&\le\log|K|+\frac{k_n}{n}\log\frac{e|H|n}{k_n},
\end{align*}
where the second inequality uses $-(1-x)\log(1-x)\le x$ for $0\le x<1$. Combining Lemma \ref{lem4c}(a) with $w=0$ and the preceding entropy bound gives
\begin{align*}
-n^{1/3} &\le n\Phi_n(v_n)\\
&=n\mr H(v_n)+n\log q_{\ell_n,H}(\mb h_n)\\
&\le k_n\log\frac{e|H|n}{k_n}-\ell_n k_n+|H|\frac{\ell_n^2k_n^2}{n}+\frac{2|H|n}{(8|H|)^{\ell_n}}.
\end{align*}

Now the inequalities $k_n > n/(|H|\ell_n^2)$, $k_n \ge 1$ and $\ell_n k_n/n<1/(4|H|)$ imply that, for all sufficiently large $n$, the right-hand side of the above inequality is bounded above by
$$
k_n\log(e|H|^2\ell_n^2)-\ell_n k_n+\frac{\ell_n k_n}{4}+o(1)
\le-\frac{\ell_n k_n}{2}+o(1)
\le-\frac12\sqrt{\frac n{|H|}}+o(1).
$$
This contradicts $n^{1/3}=o(n^{1/2})$ and proves (b). \qedhere
\end{enumerate}
\end{proof}

\begin{lem}\label{lem4f}
\begin{enumerate}
\item If $\ell_n-\log n\to\infty$, then
$$
\sum_{\mb h\in\Gen_n(H)\setminus\mc B_n}q_{\ell_n,H}(\mb h)^{n+w}\to0.
$$
\item We have
$$
\sum_{\mb h\in\Gen_n(H)\setminus\mc B_n}\bP(\mc D_n\cap\mc A_{\mb h})\to0.
$$
\end{enumerate}
\end{lem}

\begin{proof}
By Proposition \ref{prop4e}(a) with $t=n^{-2/3}$, tuples with $\Phi_n(v_{\mb h})<-n^{-2/3}$ contribute at most $(n+1)^{|H|}e^{-n^{1/3}}=o(1)$ to either sum. Every remaining tuple is covered by some $K<H$ as in Proposition \ref{prop4e}(b). Set $a_n:=n|H|e^{1-\ell_n}$. For fixed $K<H$ and $1\le k\le n/(|H|\ell_n^2)$, counting the $\binom nk(|H|-|K|)^k|K|^{n-k}$ tuples with $k$ coordinates outside $K$ and applying \eqref{eq4d} gives
$$
\sum_{\substack{\mb h\in H^n:\,k(\mb h)=k\\ \|v_{\mb h_K}-u_K\|_1\le1/(8|H|)}}q_{\ell_n,H}(\mb h)^{n+w}
< 2\binom nk\left(\frac{|H|}{|K|}-1\right)^k e^{(1-\ell_n)k}
< \frac{2a_n^k}{k!}
$$
for all sufficiently large $n$. The first inequality uses $|H|\ell_n^2k^2/n\le k$ and $(2|H|n)/(8|H|)^{\ell_n}\to0$, while the second inequality uses $\binom nk\le n^k/k!$.
Summing over $k\ge1$ and $K<H$, we have
\begin{equation} \label{eq4n2}
\sum_{\mb h\in\Gen_n(H)\setminus\mc B_n}q_{\ell_n,H}(\mb h)^{n+w}
\le o(1)+C\sum_{k\ge1}\frac{a_n^k}{k!}
=o(1)+C(e^{a_n}-1).
\end{equation}
This proves (a), since $a_n\to0$ in that case.

For (b), fix $K<H$ and $1\le k\le n/(|H|\ell_n^2)$, and consider $\mb h\in H^n$ satisfying $k(\mb h)=k$ and $\|v_{\mb h_K}-u_K\|_1\le1/(8|H|)$. Write $q:=q_{\ell_n,H}(\mb h)$. Restricting to words whose letter images all lie in $K$, \eqref{eq4b} gives
$$
q \ge \left(1-\frac{k}{n}\right)^{\ell_n} q_{\ell_n,K}(\mb h_K)
\ge\left(1-\frac{\ell_n k}{n}\right)\left(|K|^{-1}-(8|H|)^{-\ell_n}\right)
\ge \frac{1}{2|K|}
$$
for all sufficiently large $n$.

On $\mc D_n\cap\mc A_{\mb h}$, some relation contains at least two letters whose images lie outside $K$: at least one such letter occurs because $\mc D_n$ occurs and $k\ge1$, and exactly one cannot give product $1$. 
We apply the union bound over all choices of a relation and two of its letters whose images lie outside $K$. For each choice, we drop the requirement that the chosen relation maps to $1$. The independence of the relations and the inequality $q^{-1}\le2|K|$ give
$$
\bP(\mc D_n\cap\mc A_{\mb h})
\le(n+w)\binom{\ell_n}2\frac{k^2}{n^2}q^{n+w-1}
\le C\frac{\ell_n^2k^2}{n}q^{n+w}.
$$
For each fixed $K<H$ and $k$, we combine this bound with the preceding estimate $2a_n^k/k!$. Summing over $K$ and $k$ then gives
\begin{align*}
\sum_{\mb h\in\Gen_n(H)\setminus\mc B_n}\bP(\mc D_n\cap\mc A_{\mb h})
&\le o(1)+\frac{C\ell_n^2}{n}\sum_{k\ge1}\frac{k^2a_n^k}{k!}\\
&=o(1)+\frac{C\ell_n^2}{n}(a_n+a_n^2)e^{a_n}\\
&=o(1)+C\ell_n^2|H|e^{1-\ell_n}(1+a_n)e^{a_n} \to 0,
\end{align*}
which follows from $a_n=O(1)$ and $\ell_n \to \infty$. \qedhere
\end{proof}

\begin{proof}[Proof of Proposition \ref{prop4d}]
The case $H=1$ is trivial. For a nontrivial $H$, the assertion follows from \eqref{eq4a}, \eqref{eq4f}, and Lemma \ref{lem4f}(a).
\end{proof}

\subsection{Proof of Theorem \ref{thm1c}}\label{Sub43}

Let $\mc U_n\subseteq[n]$ be the set of indices of generators absent from the unreduced sampled words defining $G_{n,w}[\ell_n]$. Set $U_n:=|\mc U_n|$, so that $\mc D_n=\{U_n=0\}$. For $J\subseteq[n]$, let $\mc E_J$ be the event that $J\subseteq\mc U_n$. For $|J|=r$, we have $\bP(\mc E_J)=(1-r/n)^{\ell_n(n+w)}$. Write $(a)_r:=\prod_{i=0}^{r-1}(a-i)$ and $(a)_0:=1$. Then 
\begin{equation}\label{eq4i}
\bE((U_n)_r)=(n)_r(1-r/n)^{\ell_n(n+w)} \quad (0\le r\le n).
\end{equation}
Recall that $\lambda_s := e^{-s}$. If $\ell_n-\log n\to s\in\bR$ along a subsequence, then $\bE((U_n)_r)\to\lambda_s^r$ for every fixed integer $r\ge0$. Applying \cite[Theorem 1.22]{Bol01}, we obtain
\begin{equation}\label{eq4j}
U_n\dist Z_s\sim\Pois(\lambda_s),\quad \bP(\mc D_n)\to e^{-\lambda_s}>0.
\end{equation}

\begin{prop}\label{prop4g}
Fix $w\in\bZ_{\ge0}$ and $s\in\bR$, and let $\mc N\subseteq\bZ_{\ge1}$ be infinite. Suppose that $\ell_n$ is odd for $n\in\mc N$ and $\ell_n-\log n\to s$ as $n\to\infty$ through $\mc N$. The conditional law of $G_{n,w}[\ell_n]$, given that every generator occurs in the unreduced sampled words, converges weakly to $\mu_w$ along $\mc N$.
\end{prop}

\begin{proof}
Fix a nontrivial finite group $H$. The estimates of Section \ref{Sub42} apply along $\mc N$ with any fixed $b<s$. Fix an integer $r\ge0$ independently of $n$. Let $J\subseteq[n]$ have size $r$ and write $J^c:=[n]\setminus J$. For $\mb h\in\mc B_n$, let $\mb h_{J^c}:=(h_i)_{i\in J^c}$ with the indices in increasing order. The identity $(n-r)v_{\mb h_{J^c}}=nv_{\mb h}-\sum_{j\in J}\delta_{h_j}$ gives
$$
\|v_{\mb h_{J^c}}-v_{\mb h}\|_1
=\frac1{n-r}\left\|rv_{\mb h}-\sum_{j\in J}\delta_{h_j}\right\|_1
\le\frac{2r}{n-r}.
$$
Hence the triangle inequality gives, for all sufficiently large $n$, 
$$
\|v_{\mb h_{J^c}}-u_H\|_1
\le\|v_{\mb h}-u_H\|_1+\|v_{\mb h_{J^c}}-v_{\mb h}\|_1
\le\frac1{8|H|}+\frac{2r}{n-r}
\le\frac1{4|H|}.
$$

Conditional on $\mc E_J$, the relations are independent words on the remaining signed generators. Thus 
$$
\bP(\mc E_J\mid\mc A_{\mb h})
=\frac{\bP(\mc E_J)}{\bP(\mc A_{\mb h})}\bP(\mc A_{\mb h}\mid\mc E_J)
=\bP(\mc E_J)\left(\frac{q_{\ell_n,H}(\mb h_{J^c})}{q_{\ell_n,H}(\mb h)}\right)^{n+w}.
$$
Moreover, \eqref{eq4b} gives
$$
q_{\ell_n,H}(\mb h_{J^c})=|H|^{-1}+O((4|H|)^{-\ell_n}), \quad q_{\ell_n,H}(\mb h)=|H|^{-1}+O((8|H|)^{-\ell_n}).
$$
Combining these estimates with $(4|H|)^{-\ell_n} = o(n^{-1})$ gives, uniformly in $J$ and $\mb h$,
$$
\bP(\mc E_J)\left(\frac{q_{\ell_n,H}(\mb h_{J^c})}{q_{\ell_n,H}(\mb h)}\right)^{n+w}
=\bP(\mc E_J)(1+o(1)).
$$
Since $(U_n)_r=r!\sum_{|J|=r}\mb1_{\mc E_J}$,
$$
\bE((U_n)_r\mid\mc A_{\mb h})
=r!\sum_{\substack{J\subseteq[n]\\|J|=r}}\bP(\mc E_J\mid\mc A_{\mb h})
=(n)_r(1-r/n)^{\ell_n(n+w)}(1+o(1))
\to\lambda_s^r.
$$

By the same argument as for \eqref{eq4j}, we have $\bP(\mc D_n\mid\mc A_{\mb h})=e^{-\lambda_s}+o(1)$ uniformly in $\mb h$. By \eqref{eq4f},
\begin{equation}\label{eq4k}
\sum_{\mb h\in\mc B_n}\bP(\mc D_n\cap\mc A_{\mb h})
=(e^{-\lambda_s}+o(1))\sum_{\mb h\in\mc B_n}q_{\ell_n,H}(\mb h)^{n+w}
\to e^{-\lambda_s}|H|^{-w}.
\end{equation}
Lemma \ref{lem4f}(b), \eqref{eq4j} and \eqref{eq4k} imply that
$$
\bE(\#\Sur(G_{n,w}[\ell_n],H)\mid\mc D_n)
=\frac{\sum_{\mb h\in\Gen_n(H)}\bP(\mc D_n\cap\mc A_{\mb h})}{\bP(\mc D_n)}
=\frac{e^{-\lambda_s}|H|^{-w}+o(1)}{e^{-\lambda_s}+o(1)}
\to|H|^{-w}
$$
for every nontrivial $H$. The case $H=1$ is trivial, so Propositions \ref{prop2a} and \ref{prop2b} finish the proof.
\end{proof}

\begin{proof}[Proof of Theorem \ref{thm1c}]
Take $w=u$.
\begin{enumerate}
\item If $\ell_n$ is eventually odd and $s_n\to\infty$, Proposition \ref{prop4d} gives the convergence of moments. Propositions \ref{prop2a} and \ref{prop2b} give $G_{n,u}\dist\mu_u$, proving the if part of (a).

\item Work along the subsequence in (b). By \eqref{eq4j}, $U_n\dist Z_s$. Conditional on $U_n=r$, we may relabel the unused generators as $x_{n-r+1},\ldots,x_n$. These generators occur in no relation and hence contribute a free profinite factor $\Fhat_r$. Thus
$$
\mc L(G_{n,u}\mid U_n=r)=\mc L(\Fhat_r\amalg V_{n,r}),
$$
where $V_{n,r}$ has the law of $G_{n-r,u+r}[\ell_n]$ conditioned on every generator occurring, since the remaining presentation has $n-r$ generators and $n+u=(n-r)+(u+r)$ relations. For fixed $r$, $\ell_n-\log(n-r)\to s$, so Proposition \ref{prop4g} gives $V_{n,r}\dist X_{u+r}$.

Fix a finite set $\mc S$ of isomorphism classes of finite groups. Every continuous homomorphism $G\to B$ with $B\in\mc S^{\mr{cl}}$ factors uniquely through $G^{\mc S}$. By the universal property of the free profinite product \cite[Section 9.1]{RZ10}, every continuous homomorphism $\Fhat_r\amalg G\to B$ factors uniquely through the natural quotient $\Fhat_r\amalg G\to\Fhat_r\amalg G^{\mc S}$. This holds for every $B\in\mc S^{\mr{cl}}$, so the definition of the $\mc S$-completion gives
$$
(\Fhat_r\amalg G)^{\mc S}\cong(\Fhat_r\amalg G^{\mc S})^{\mc S}.
$$
If $G\in\mc P$, then $G^{\mc S}$ is finite and $\Fhat_r\amalg G^{\mc S}$ is finitely generated. As recalled in Section \ref{Sub22}, every finitely generated profinite group belongs to $\mc P$, so the right-hand side is finite. Since this holds for every finite $\mc S$, we obtain $\Fhat_r\amalg G\in\mc P$.

To prove the continuity of the map $G\mapsto\Fhat_r\amalg G$, it suffices to show that the inverse image of each $U_{\mc S,A}$ is open. For every finite $\mc S$ and every $A\in\mc S^{\mr{cl}}$, the displayed isomorphism gives
$$
\{G\in\mc P:(\Fhat_r\amalg G)^{\mc S}\cong A\}
=\bigcup_{\substack{B\in\mc S^{\mr{cl}}\\(\Fhat_r\amalg B)^{\mc S}\cong A}}U_{\mc S,B}.
$$
The left-hand side is the inverse image of $U_{\mc S,A}$, and the right-hand side is a union of open sets. Thus the map $G\mapsto\Fhat_r\amalg G$ is continuous. Applying the continuous mapping theorem to $V_{n,r}\dist X_{u+r}$ therefore gives $\Fhat_r\amalg V_{n,r}\dist\Fhat_r\amalg X_{u+r}$ for each fixed $r$.
For every bounded continuous function $f:\mc P\to\bR$ and fixed integer $L\ge0$, the convergence for each fixed $r$ gives
$$
\sum_{r=0}^{L}\bP(U_n=r)\bE(f(G_{n,u})\mid U_n=r) \to
\sum_{r=0}^{L}\bP(Z_s=r)\bE(f(\Fhat_r\amalg X_{u+r})).
$$
Since $U_n\dist Z_s$, we have
\begin{align*}
\limsup_{n\to\infty}\left|\bE(f(G_{n,u}))-\bE(f(\mc Y_{u,s}))\right| 
&\le\|f\|_\infty\limsup_{n\to\infty}\left(\bP(U_n>L)+\bP(Z_s>L)\right) \\
&=2\|f\|_\infty\bP(Z_s>L).
\end{align*}
Taking $L\to\infty$ proves (b).

\item If $s_n\to-\infty$, then $\ell_n\le\log n$ for all sufficiently large $n$. Applying \eqref{eq4i} gives
\begin{align*}
\log\bE(U_n)=\log n+\ell_n(n+u)\log(1-1/n)=-s_n+o(1),\\
\frac{\bE((U_n)_2)}{(\bE(U_n))^2}=\left(1-\frac1n\right)
\left(\frac{1-2/n}{(1-1/n)^2}\right)^{\ell_n(n+u)}\le1.
\end{align*}

Thus $\Var(U_n)\le\bE(U_n)$ and $\bE(U_n)\to\infty$, so Chebyshev's inequality gives $U_n\xrightarrow{\bP}\infty$. Sending every generator used in the unreduced relators defining $G_{n,u}$ to $0\in C_2$ and assigning arbitrary images to the unused generators gives
$$
\#\Sur(G_{n,u},C_2)\ge2^{U_n}-1\xrightarrow{\bP}\infty.
$$
Since $G\mapsto\#\Sur(G,C_2)$ is locally constant on $\mc P$, a weakly convergent subsequence of the group laws would give a tight subsequence of these counts, a contradiction.

\item If $\ell_n$ is even, sending every generator to $1\in C_2$ sends each relation to $\ell_n\bmod2=0$ and gives a surjection $G_{n,u}\to C_2$.
\end{enumerate}

We now prove the only if part of (a). For the Haar-relator quotient $Q_{N,u}$, the event $\#\Sur(Q_{N,u},C_2)=0$ is equivalent to $N+u$ independent uniform vectors spanning $C_2^N$. Since the event $\#\Sur(G,C_2)=0$ is clopen on $\mc P$, the convergence $Q_{N,u}\dist X_u$ in Proposition \ref{prop2b} gives
\begin{align*}
a_u:=\bP(\#\Sur(X_u,C_2)=0)
&=\lim_{N\to\infty}\bP(\#\Sur(Q_{N,u},C_2)=0)\\
&=\lim_{N\to\infty}\prod_{j=u+1}^{N+u}(1-2^{-j})\\
&=\prod_{j=u+1}^{\infty}(1-2^{-j})>0.
\end{align*}
Thus convergence to $\mu_u$ forces $\ell_n$ to be eventually odd. Moreover, $\Fhat_r\amalg G$ surjects onto $C_2$ whenever $r\ge1$, so
$$
\bP(\#\Sur(\mc Y_{u,s},C_2)=0)=e^{-\lambda_s}a_u<a_u.
$$
Hence $\mc L(\mc Y_{u,s})\ne\mu_u$. If $s_n\not\to\infty$, pass to a subsequence with $s_n\to s\in\bR$ or $s_n\to-\infty$. Parts (b) and (c) then rule out the convergence to $\mu_u$, proving the only if part of (a).
\end{proof}

\section{Words with balanced integer exponents}\label{Sec5}

Let $\mb h\in H^n$. Take independent integer-valued random variables $\xi_1,\ldots,\xi_n$, each $\alpha$-balanced modulo every prime dividing $|H|$, and an independent uniform permutation $\sigma\in S_n$. Set
$$
q^{\mr{ord}}_H(\mb h)
:=\bP(h_{\sigma(1)}^{\xi_{\sigma(1)}}\cdots
h_{\sigma(n)}^{\xi_{\sigma(n)}}=1).
$$
The estimates for $q_H^{\mr{ord}}$ are uniform in the laws of the $\xi_i$ and in $0<\alpha\le1/2$.

\begin{prop}\label{prop5a}
Fix $b>0$. There exist a constant $c>0$ and a positive integer $n_0$, depending only on $H$ and $b$, with the following property. Let $n\ge n_0$, and suppose that $K\le H$ has a generating set $S$ each of whose elements occurs at least $bn$ times among the coordinates of $\mb h \in H^n$. Then, with $k=k(\mb h):=\#\{i:h_i\notin K\}$,
\begin{equation}\label{eq5a}
q^{\mr{ord}}_H(\mb h)
\le\frac1{|K|}-\frac{1-e^{-c\alpha k}}{4|K|}+e^{-c\alpha n}.
\end{equation}
If $k=0$, then in addition
\begin{equation}\label{eq5b}
\left|q^{\mr{ord}}_H(\mb h)-\frac1{|K|}\right| \le e^{-c\alpha n}.
\end{equation}
\end{prop}

\subsection{Return probability estimates}\label{Sub51}

For $K\le H$, the operator $\Pi_K:=T_{u_K}$ is the orthogonal projection onto the right $K$-invariant functions. We write $\Pi_g:=\Pi_{\langle g\rangle}$ and $x_+:=\max\{x,0\}$. For a probability vector $f$ on $H$, write $f(A):=\sum_{g\in A}f(g)$ for $A\subseteq H$. Then $(\Pi_Kf)(x)=f(xK)/|K|$, so
\begin{equation} \label{eq5n1}
\|\Pi_Kf\|_2^2=\frac1{|K|}\sum_{xK\in H/K}f(xK)^2\le\frac1{|K|}.
\end{equation}

\begin{lem}\label{lem5b}
Let $g\in H$, and let $\xi$ be an integer-valued random variable that is $\alpha$-balanced modulo every prime dividing $|H|$. For $\nu:=\mc L(g^\xi)$,
\begin{equation}\label{eq5c}
\|f\|_2^2-\|T_\nu f\|_2^2
\ge\frac{\alpha}{2|H|^3}\|f-\Pi_gf\|_2^2
\quad\text{for all }f\in\ell^2(H).
\end{equation}
The same bound holds for $\check\nu$.
\end{lem}

\begin{proof}
If $g=1$, then both sides of \eqref{eq5c} vanish. Otherwise, set $L:=\langle g\rangle$ and $d:=|L|>1$. A coset of a proper subgroup of $L$ corresponds to an event $\xi\equiv a\pmod r$, where $r>1$ divides $d$. For a prime $p\mid r$, we have
$$
\bP(\xi\equiv a\pmod r)\le\bP(\xi\equiv a\pmod p)\le1-\alpha.
$$
Thus $\nu$ is $\alpha$-coset-balanced on $L$. For each left coset $xL$, the function $a\mapsto(f-\Pi_gf)(xa)$ belongs to $\ell^2_0(L)$. Applying Proposition \ref{prop2c} to these functions and summing the squared bounds over the cosets gives
$$
\|T_\nu(f-\Pi_gf)\|_2\le e^{-\alpha/(2d^3)}\|f-\Pi_gf\|_2.
$$
Since $T_\nu$ fixes each element of the image of $\Pi_g$ and preserves its orthogonal complement, the orthogonal decomposition $f=\Pi_gf+(f-\Pi_gf)$ gives
\begin{align*}
\|f\|_2^2-\|T_\nu f\|_2^2
&=\|f-\Pi_gf\|_2^2-\|T_\nu(f-\Pi_gf)\|_2^2\\
&\ge(1-e^{-\alpha/d^3})\|f-\Pi_gf\|_2^2\\
&\ge\frac{\alpha}{2|H|^3}\|f-\Pi_gf\|_2^2,
\end{align*}
where the last inequality uses $d\le|H|$ and $1-e^{-x}\ge x/2$ for $0\le x\le1$. This proves \eqref{eq5c}. Replacing $\xi$ by $-\xi$ gives the assertion for $\check\nu$.
\end{proof}

\begin{lem}\label{lem5c}
There exists $\gamma_H>0$ such that, for every $K\le H$, the following hold.
\begin{enumerate}
\item If $S\subseteq K$ generates $K$, then
$$
\sum_{a\in S}\|f-\Pi_af\|_2^2
\ge\gamma_H\|f-\Pi_Kf\|_2^2
\quad\text{for all }f\in\ell^2(H).
$$
\item If $g\notin K$ and $f$ is a probability vector on $H$, then
$$
\|f-\Pi_Kf\|_2^2+\|f-\Pi_gf\|_2^2
\ge\gamma_H\left(\|f\|_2^2-\frac1{2|K|}\right)_+.
$$
\end{enumerate}
\end{lem}

\begin{proof}
Part (a) is trivial when $K=1$. For a nontrivial $K$, $\left\langle (I-\Pi_a)f, \Pi_a f \right\rangle=0$ for each $a \in S$, so
$$
\left\langle\sum_{a\in S}(I-\Pi_a)f,f\right\rangle
=\sum_{a\in S}\|f-\Pi_af\|_2^2.
$$
Consider the positive semidefinite operator $A := \sum_{a\in S}(I-\Pi_a)$. Since $S$ generates $K$, we have
\begin{align*}
f\in\ker A 
& \quad\Leftrightarrow\quad \langle Af,f\rangle=0\\
& \quad\Leftrightarrow\quad \Pi_a f = f \quad\text{for every } a \in S \\
& \quad\Leftrightarrow\quad f\text{ is right }K\text{-invariant}\\
& \quad\Leftrightarrow\quad f\in\mr{im}\Pi_K.
\end{align*}
Thus $A$ is positive definite on the orthogonal complement of $\mr{im}\Pi_K$. Writing $\lambda_{K,S}>0$ for its least positive eigenvalue, we obtain
$$
\sum_{a\in S}\|f-\Pi_af\|_2^2
=\langle A(f-\Pi_Kf),f-\Pi_Kf\rangle
\ge\lambda_{K,S}\|f-\Pi_Kf\|_2^2.
$$
Taking the minimum of $1$ and the finitely many $\lambda_{K,S}$ gives a common $\gamma_H>0$ in (a).

For (b), let $L:=\langle K,g\rangle$, the subgroup generated by $K\cup\{g\}$. For $a\in K$, the inclusion $\mr{im}\Pi_K\subseteq\mr{im}\Pi_a$ gives $\|f-\Pi_af\|_2\le\|f-\Pi_Kf\|_2$. Applying (a) to the generating set $K\cup\{g\}$ of $L$, we obtain
$$
\|f-\Pi_Kf\|_2^2+\|f-\Pi_gf\|_2^2
\ge\frac1{|K|}\sum_{a\in K\cup\{g\}}\|f-\Pi_af\|_2^2
\ge\frac{\gamma_H}{|K|}\|f-\Pi_Lf\|_2^2.
$$
Since $g\notin K$, \eqref{eq5n1} gives $\|\Pi_Lf\|_2^2\le|L|^{-1}\le(2|K|)^{-1}$. As $\Pi_L$ is an orthogonal projection,
$$
\|f-\Pi_Lf\|_2^2=\|f\|_2^2-\|\Pi_Lf\|_2^2 \ge\left(\|f\|_2^2-\frac1{2|K|}\right)_+.
$$
Since $|K|\le|H|$, replacing $\gamma_H$ by $\gamma_H/|H|$ proves (b) while preserving (a).
\end{proof}

\begin{proof}[Proof of Proposition \ref{prop5a}]
Set $\nu_i:=\mc L(h_i^{\xi_i})$, $f_0:=\delta_1$, and
$$
f_t:=\nu_{\sigma(1)}*\cdots*\nu_{\sigma(t)},\quad F_t:=\|f_t\|_2^2.
$$
Conditional on $\sigma$, $f_t$ is the law of the product of the first $t$ factors. In particular, $q_H^{\mr{ord}}(\mb h)=\bE(f_n(1))$. For $0\le t\le n$, let $\mc I_t:=[n]\setminus\{\sigma(1),\ldots,\sigma(t)\}$ be the set of remaining indices. For $0 \le t< n$, we have $f_{t+1}=f_t*\nu_{\sigma(t+1)}=T_{\check\nu_{\sigma(t+1)}}f_t$ and $0\le F_{t+1}\le F_t\le1$. Set
$$
\tau_0:=\min\{b/4,1/(8|H|)\},\quad m:=\lfloor\tau_0n\rfloor.
$$
Increasing $n_0$ if necessary, we may assume $m\ge\tau_0n/2$.

\emph{Step 1: reduction to an energy bound.}
Set $r:=\lfloor n/2\rfloor$ and
$$
\nu_-:=\nu_{\sigma(1)}*\cdots*\nu_{\sigma(r)},\quad
\nu_+:=\nu_{\sigma(r+1)}*\cdots*\nu_{\sigma(n)}.
$$Each contains at least $2m$ factors, and its first $2m$ indices have the same marginal law as $(\sigma(1),\ldots,\sigma(2m))$. Convolution with a probability measure does not increase the $\ell^2$-norm, so 
\begin{align*}
\bE(\|\nu_-\|_2^2)
&\le\bE(\|\nu_{\sigma(1)}*\cdots*\nu_{\sigma(2m)}\|_2^2)
=\bE(F_{2m}),\\
\bE(\|\nu_+\|_2^2)
&\le\bE(\|\nu_{\sigma(r+1)}*\cdots*\nu_{\sigma(r+2m)}\|_2^2)
=\bE(F_{2m}).
\end{align*}
Hence the Cauchy--Schwarz inequality gives
\begin{equation} \label{eq5n2}
q_H^{\mr{ord}}(\mb h)
=\bE((\nu_- * \nu_+)(1))
=\bE(\langle\nu_-,\check\nu_+\rangle)
\le(\bE(\|\nu_-\|_2^2))^{1/2}(\bE(\|\nu_+\|_2^2))^{1/2} 
\le\bE(F_{2m}).
\end{equation}

\emph{Step 2: the one-step energy decrease.}
For each $t<n$, fix the values of $\sigma(1),\ldots,\sigma(t)$. Then $f_t$ is determined, and $\sigma(t+1)$ is uniform on $\mc I_t$. Averaging Lemma \ref{lem5b} over $\mc I_t$ gives
$$
\bE(F_t-F_{t+1}\mid \sigma(1),\ldots,\sigma(t))
\ge\frac{\alpha}{2|H|^3(n-t)}\Delta_t,
\quad \Delta_t:=\sum_{i\in\mc I_t}\|f_t-\Pi_{h_i}f_t\|_2^2.
$$
For $t<2m$, each $a\in S$ still occurs at least $bn/2$ times, so Lemma \ref{lem5c}(a) gives
$$
\Delta_t\ge\frac{b\gamma_Hn}{2}\|f_t-\Pi_Kf_t\|_2^2
+\sum_{\substack{i\in\mc I_t\\h_i\notin K}}
\|f_t-\Pi_{h_i}f_t\|_2^2.
$$
Set $k_t:=\#\{i\in\mc I_t:h_i\notin K\}$. Since $k_t\le n$, summing Lemma \ref{lem5c}(b) over the indices $i \in \mc I_t$ with $h_i \notin K$ gives
$$
\Delta_t\ge\gamma_H\min\{b\gamma_H/2,1\}k_t
\left(F_t-\frac1{2|K|}\right)_+.
$$
Since $\Pi_K$ is an orthogonal projection, \eqref{eq5n1} gives
$$
\|f_t-\Pi_Kf_t\|_2^2=F_t-\|\Pi_Kf_t\|_2^2\ge F_t-|K|^{-1}.
$$
Combining this with the preceding bounds gives constants $a_0,a_1\in(0,1]$, depending only on $H$ and $b$, such that for each $t<2m$,
\begin{equation}\label{eq5f}
\begin{aligned}
\bE(F_t-F_{t+1}\mid\sigma(1),\ldots,\sigma(t)) &\ge a_0\alpha(F_t-|K|^{-1}),\\
\bE(F_t-F_{t+1}\mid\sigma(1),\ldots,\sigma(t)) &\ge a_1\alpha\frac{k_t}{n}(F_t-(2|K|)^{-1})_+.
\end{aligned}
\end{equation}

\emph{Step 3: the first energy bound.}
The first inequality of \eqref{eq5f} gives, for $0\le t<m$,
$$
\bE(F_{t+1})-|K|^{-1} \le(1-a_0\alpha)(\bE(F_t)-|K|^{-1}).
$$
Iterating from $F_0=1$ gives, for some constant $c>0$,
\begin{equation}\label{eq5g}
\bE(F_m)
\le |K|^{-1}+(1-|K|^{-1})(1-a_0\alpha)^m
\le |K|^{-1}+e^{-c\alpha n}.
\end{equation}
If $k=0$, then $f_n$ is supported on $K$, so $\|f_n-u_K\|_2^2=F_n-|K|^{-1}$. Since $F_n\le F_m$, the Cauchy--Schwarz inequality and \eqref{eq5g} give
$$
\left|q_H^{\mr{ord}}(\mb h)-|K|^{-1}\right|
\le(\bE(\|f_n-u_K\|_2^2))^{1/2}
\le(\bE(F_m)-|K|^{-1})^{1/2}
\le e^{-c\alpha n/2}.
$$
This proves \eqref{eq5b} after replacing $c$ with $c/2$.

\emph{Step 4: the additional loss when $k\ge1$.}
Set $X_t:=(F_t-(2|K|)^{-1})_+$. Each of the $k$ indices with $h_i\notin K$ occurs among the first $t$ entries of $\sigma$ with probability $t/n$, so $\bE(k-k_t)=kt/n$. For $m\le t<2m$, the bounds $0\le X_t\le1$ and $t/n\le1/(4|K|)$ give
$$
\bE(k_tX_t)
\ge k\bE(X_t)-\bE(k-k_t)
\ge k\left(\bE(F_t)-\frac1{2|K|}-\frac tn\right)
\ge k\left(\bE(F_t)-\frac3{4|K|}\right).
$$
The second inequality of \eqref{eq5f} gives, for $m\le t<2m$,
$$
\bE(F_{t+1})-\frac3{4|K|}
\le(1-a_1\alpha k/n)\left(\bE(F_t)-\frac3{4|K|}\right).
$$
Iterating from \eqref{eq5g} gives, for some constant $c>0$ depending only on $H$ and $b$,
$$
\bE(F_{2m}) \le\frac3{4|K|}+\left(\frac1{4|K|}+e^{-c\alpha n}\right)
(1-a_1\alpha k/n)^m
\le\frac1{|K|}-\frac{1-e^{-c\alpha k}}{4|K|}+e^{-c\alpha n}.
$$
Combining the above inequality with \eqref{eq5n2} proves \eqref{eq5a} for $k \ge 1$.
\end{proof}

\subsection{Proof of Theorem \ref{thm1d}}\label{Sub53}

We may assume that $H$ is nontrivial. Set $\alpha_n:=\min_{p\mid|H|}\alpha_{n,p}$, so that $n\alpha_n/\log n\to\infty$. For $j\in[n+u]$, write $q^{\mr{ord}}_{n,j,H}(\mb h):=\bP(\varphi_{\mb h}(r^{\mr{ord}}_{n,j})=1)$. Applying \eqref{eq2n1} gives
\begin{equation}\label{eq5h}
\bE(\#\Sur(G^{\mr{ord}}_{n,u},H))
=\sum_{\mb h\in\Gen_n(H)}
\prod_{j=1}^{n+u}q^{\mr{ord}}_{n,j,H}(\mb h).
\end{equation}
For each $j$, we apply Proposition \ref{prop5a} with $\xi_i=\xi_{i,j}$ and $\alpha=\alpha_n$. As in Section \ref{Sub41}, write $k(\mb h):=\#\{i:h_i\notin K\}$ for a fixed $K\le H$.
Let $b:=(2|H|)^{-1}$, $c>0$ be a constant given by Proposition \ref{prop5a} for this $b$, and $d:=(1-e^{-c})/4$. Let $\mr h(x):=-x\log x-(1-x)\log(1-x)$ be the binary entropy function for $0\le x\le1$. Choose $0<\delta<1/2$ such that
\begin{equation}\label{eq5i}
\mr h(\delta)+\delta\log(|H|-1)<\frac d2.
\end{equation}

For each $K\le H$, let $\mc N_K$ be the open $\ell^1$-ball of radius $\min\{\delta,b\}$ centered at $u_K$ in the space of probability vectors on $H$. These neighborhoods are pairwise disjoint, since
\begin{align*}
\|u_K-u_L\|_1 &= \sum_{x \in H}|u_K(x)-u_L(x)| \\
&= \sum_{x \in H} (u_K(x)+u_L(x)-2\min\{u_K(x),u_L(x)\}) \\
&= 2-2\frac{|K \cap L|}{\max\{|K|,|L|\}} \ge 1
\end{align*}
for distinct subgroups $K,L\le H$ and $b\le1/4$. Every $v\in\mc N_K$ satisfies
\begin{equation}\label{eq5j}
v(a)\ge b\ \text{ for all }a\in K,
\quad
\sum_{a\notin K}v(a)\le \|v-u_K\|_1 < \delta.
\end{equation}

The following cases define a partition of $\Gen_n(H)$. In the last two cases, $K$ ranges over the proper subgroups of $H$.
\begin{itemize}[leftmargin=3.2em,itemsep=2pt,topsep=3pt,parsep=0pt]
\item[$(\mc D_1)$] $v_{\mb h}\notin\bigcup_{K\le H}\mc N_K$;
\item[$(\mc D_2)$] $v_{\mb h}\in\mc N_H$;
\item[$(\mc D_3^K)$] $v_{\mb h}\in\mc N_K$ and $\alpha_nk(\mb h)\le1$;
\item[$(\mc D_4^K)$] $v_{\mb h}\in\mc N_K$ and $\alpha_nk(\mb h)>1$.
\end{itemize}
Note that $k(\mb h)\ge1$ in cases $(\mc D_3^K)$ and $(\mc D_4^K)$, since $\mb h$ generates $H$. For a condition $\mc D$, define
$$
G_{\mc D}(n):=\sum_{\substack{\mb h\in\Gen_n(H)\\ \mc D\text{ is satisfied}}}
\prod_{j=1}^{n+u}q^{\mr{ord}}_{n,j,H}(\mb h).
$$
Then \eqref{eq5h} becomes
$$
\bE(\#\Sur(G^{\mr{ord}}_{n,u},H))
=G_{\mc D_1}(n)+G_{\mc D_2}(n)
+\sum_{K<H}(G_{\mc D_3^K}(n)+G_{\mc D_4^K}(n)).
$$
We prove that $G_{\mc D_2}(n)$ gives the main term $|H|^{-u}$, while all remaining terms tend to zero.

\begin{prop}\label{prop5d1}
$\lim_{n\to\infty}G_{\mc D_1}(n)=0$.
\end{prop}

\begin{proof}
For each type $v$ of a tuple in $H^n$, set
$$
Q_n(v):=\max\{q^{\mr{ord}}_{n,j,H}(\mb h):j\in[n+u],\ v_{\mb h}=v\}.
$$
If $v_n\to v$, set $S:=\supp(v)$ and $b_v:=\frac12\min_{a\in S}v(a)>0$. Then $v_n(a)\ge b_v$ for every $a\in S$ and all sufficiently large $n$. Proposition \ref{prop5a}, applied with frequency parameter $b_v$ and $K=\langle S\rangle$, gives
$$
Q_n(v_n)\le|K|^{-1}+e^{-c_v\alpha_n n},
$$
where $c_v>0$ depends only on $H$ and $v$. The bound is uniform in $j$ and in tuples of type $v_n$. Since $\alpha_n n\to\infty$,
$$
\limsup_{n\to\infty}Q_n(v_n)\le|\langle\supp(v)\rangle|^{-1}.
$$
By Lemma \ref{lem2d}, there exists $c>0$ such that, for all sufficiently large $n$,
$$
\mr H(v)+\log Q_n(v)\le-c
\quad\text{for every type }v\notin\bigcup_{K\le H}\mc N_K.
$$
Condition $(\mc D_1)$ ensures that this bound applies to $v_{\mb h}$. Since $\prod_jq^{\mr{ord}}_{n,j,H}(\mb h)\le Q_n(v_{\mb h})^n$, each type $v$ is realized by at most $e^{n\mr H(v)}$ tuples, and there are at most $(n+1)^{|H|}$ types, we obtain
\begin{equation*}
G_{\mc D_1}(n)
\le\sum_{\substack{v\text{ a type of length }n\\v\notin\bigcup_{K\le H}\mc N_K}} e^{n\mr H(v)}Q_n(v)^n
\le(n+1)^{|H|}e^{-cn}\to0. \qedhere
\end{equation*}
\end{proof}

\begin{prop}\label{prop5d2}
$\lim_{n\to\infty}G_{\mc D_2}(n)=|H|^{-u}$.
\end{prop}

\begin{proof}
If $v_{\mb h}\in\mc N_H$, every element of $H$ occurs at least $bn$ times in $\mb h$. Then \eqref{eq5b} gives
$$
q^{\mr{ord}}_{n,j,H}(\mb h)=|H|^{-1}(1+O(e^{-c\alpha_n n}))
$$
uniformly in $j$ and $\mb h$. Since $n\alpha_n/\log n\to\infty$,
$$
\prod_{j=1}^{n+u}q^{\mr{ord}}_{n,j,H}(\mb h)
=|H|^{-n-u}\exp(O(ne^{-c\alpha_n n}))
=|H|^{-n-u}(1+o(1)).
$$
By \eqref{eq2b}, there are $|H|^n(1+o(1))$ such tuples, all generating $H$. Hence $G_{\mc D_2}(n)=|H|^{-u}(1+o(1))$.
\end{proof}

\begin{prop}\label{prop5e}
$\lim_{n\to\infty}G_{\mc D_3^K}(n)=\lim_{n\to\infty}G_{\mc D_4^K}(n)=0$ for every $K<H$.
\end{prop}

\begin{proof}
Fix $K<H$, set $D:=|H|/|K|-1\ge1$, and let $\mb h$ satisfy $(\mc D_3^K)$ or $(\mc D_4^K)$, with $k:=k(\mb h)$. By \eqref{eq5j}, every element of $K$ occurs at least $bn$ times in $\mb h$, so Proposition \ref{prop5a} with $S=K$, together with the inequality $1-e^{-cy}\ge(1-e^{-c})\min\{y,1\}$ for $y\ge0$, gives
$$
q^{\mr{ord}}_{n,j,H}(\mb h)
\le \frac{1}{|K|}-\frac{d}{|K|}\min\{\alpha_n k,1\}+e^{-c\alpha_n n}
\quad\text{for all } j\in[n+u].
$$
Since $1+x\le e^x$ and $|K|(n+u)e^{-c\alpha_n n}\to0$, we obtain for all sufficiently large $n$,
\begin{align*}
\prod_{j=1}^{n+u}q^{\mr{ord}}_{n,j,H}(\mb h)
&\le |K|^{-(n+u)}\left(1-d\min\{\alpha_n k,1\}+|K|e^{-c\alpha_n n}\right)^{n+u}\\
&\le |K|^{-(n+u)}\exp\left(-d(n+u)\min\{\alpha_n k,1\}+|K|(n+u)e^{-c\alpha_n n}\right)\\
&\le 2|K|^{-(n+u)}e^{-dn\min\{\alpha_n k,1\}},
\end{align*}
uniformly in $\mb h$. At most $\binom nk(|H|-|K|)^k|K|^{n-k}=\binom nkD^k|K|^n$ tuples have exactly $k$ coordinates outside $K$. 

In case $(\mc D_3^K)$, we have $\min\{\alpha_n k,1\}=\alpha_n k$. By the binomial theorem, summing over all $1\le k\le n$ gives
\begin{equation}\label{eq5k}
\begin{aligned}
G_{\mc D_3^K}(n)
&\le2|K|^{-u}\left((1+De^{-d\alpha_n n})^n-1\right)\\
&\le2|K|^{-u}\left(\exp(nDe^{-d\alpha_n n})-1\right)\to0,
\end{aligned}
\end{equation}
since $nDe^{-d\alpha_n n}=De^{\log n-d\alpha_n n}\to0$.
In case $(\mc D_4^K)$, we have $\min\{\alpha_n k,1\}=1$ so
\begin{equation}\label{eq5l}
G_{\mc D_4^K}(n)\le2|K|^{-u}e^{-dn}\sum_{0\le k\le\delta n}\binom nkD^k.
\end{equation}

For \eqref{eq5l}, the function $y\mapsto \mr h(y)+y\log D$ is increasing on $[0,\delta]$. Using the binary entropy bound $\binom nk\le e^{n\mr h(k/n)}$, we obtain
$$
\sum_{0\le k\le\delta n}\binom nkD^k
\le(n+1)\exp\left(n\sup_{0\le y\le\delta}(\mr h(y)+y\log D)\right)
=(n+1)e^{n(\mr h(\delta)+\delta\log D)}.
$$
Since $D\le|H|-1$, \eqref{eq5i} gives $G_{\mc D_4^K}(n)\le2|K|^{-u}(n+1)e^{-dn/2}\to0$.
\end{proof}

Propositions \ref{prop5d1}, \ref{prop5d2}, and \ref{prop5e} give $\bE(\#\Sur(G^{\mr{ord}}_{n,u},H)) \to |H|^{-u}$. By Propositions \ref{prop2a} and \ref{prop2b}, we have $G^{\mr{ord}}_{n,u}\dist\mu_u$, which completes the proof of Theorem \ref{thm1d}.

\section{Abelianization and integral cokernels}\label{Sec6}

For a finitely generated abelian group $B$ and a finite set $P$ of primes, write $B(P):=\prod_{p\in P}(B\otimes_{\bZ}\bZ_p)$. The following corollary recovers the finite-prime cokernel law of Nguyen and Wood \cite[Theorem 4.1]{NW22} from Theorem \ref{thm1d} by abelianization.

\begin{cor}\label{cor6a}
Under the hypotheses of Theorem \ref{thm1d},
$$
\lim_{n\to\infty}\bP(B_n(P)\cong A)
=\frac1{|A|^u|\Aut(A)|}
\prod_{p\in P}\prod_{k=u+1}^{\infty}(1-p^{-k})
$$
for every finite abelian group $A$ and every finite set $P$ of primes including all those dividing $|A|$.
\end{cor}

\begin{proof}
The case $P=\varnothing$ is immediate, so assume that $P$ is nonempty. The image of $r^{\mr{ord}}_{n,j}$ in $\widehat{\bZ}^n$ is the $j$-th column of $\Xi_n$. Since $\Xi_n\widehat{\bZ}^{n+u}$ is compact in $\widehat{\bZ}^n$ and $\bZ^{n+u}$ is dense in $\widehat{\bZ}^{n+u}$, we have $\overline{\Xi_n\bZ^{n+u}}=\Xi_n\widehat{\bZ}^{n+u}$. Here $\overline{\Xi_n\bZ^{n+u}}$ denotes the closure of $\Xi_n\bZ^{n+u}$ in $\widehat{\bZ}^n$. Thus abelianization and the Smith normal form of $\Xi_n$ give
$$
(G^{\mr{ord}}_{n,u})^{\mr{ab}}
\cong\widehat{\bZ}^n/\overline{\Xi_n\bZ^{n+u}}
=\widehat{\bZ}^n/\Xi_n\widehat{\bZ}^{n+u}
\cong\widehat{B_n}.
$$
Write $A=\prod_{p\in P}A_p$ ($A_p$ denotes the $p$-primary component of $A$), choose $e_p\ge0$ with $p^{e_p}A_p=0$, and set $m:=\prod_{p\in P}p^{e_p+1}$. Following the argument of \cite[Corollary 3.4]{Woo19}, we deduce that
$$
B(P)\cong A \quad\Leftrightarrow\quad B/mB\cong A
$$
for every finitely generated abelian group $B$. 

The closure $\{C_m\}^{\mr{cl}}$ consists of the finite abelian groups annihilated by $m$. Consequently, the pro-$\{C_m\}^{\mr{cl}}$ completion of a finitely generated profinite group $G$ is its maximal abelian quotient annihilated by $m$, namely $G^{\mr{ab}}/mG^{\mr{ab}}$. Applying this to $G^{\mr{ord}}_{n,u}$ gives
$$
(G^{\mr{ord}}_{n,u})^{\{C_m\}} \cong \widehat{B_n}/m\widehat{B_n} \cong B_n/mB_n.
$$
The preceding equivalence and the definition of $U_{\mc S,A}$ imply that
$$
B_n(P)\cong A
\quad\Leftrightarrow\quad (G^{\mr{ord}}_{n,u})^{\{C_m\}}\cong A
\quad\Leftrightarrow\quad G^{\mr{ord}}_{n,u}\in U_{\{C_m\},A}.
$$
Since $U_{\{C_m\},A}$ is clopen in $\mathcal{P}$, Theorem \ref{thm1d} gives
$$
\lim_{n\to\infty}\bP(B_n(P)\cong A)=\mu_u(U_{\{C_m\},A}).
$$

For each $p\in P$, we have $(m/p)A=0$, so every abelian extension of $A$ by $C_p$ is annihilated by $m$. Thus \cite[(3.2)]{LW20}, with the same calculation as in \cite[Example 12.3]{LW20}, gives
\begin{equation*}
\mu_u(U_{\{C_m\},A}) =\frac1{|A|^u|\Aut(A)|}\prod_{p\in P}\prod_{k=u+1}^{\infty}(1-p^{-k}). \qedhere
\end{equation*}
\end{proof}

\section*{Acknowledgments}
Jiwan Jung was supported by the National Research Foundation of Korea (NRF) grant funded by the Korea government (MSIT) (No. RS-2026-25508638). Jungin Lee was supported by the National Research Foundation of Korea (NRF) grants funded by the Korea government (MSIT) (No. RS-2024-00334558 and No. RS-2025-02262988).

\section*{Statement on AI use}

OpenAI’s ChatGPT 5.6 Sol and 6 Astra generated the proofs of the main results. The models also proposed the statement of Theorem \ref{thm1b} and an initial version of Theorem \ref{thm1c} with Poisson-distributed word lengths. After discussion with Anthropic’s Claude Fable 5.1, the authors reformulated Theorem \ref{thm1c} for words of fixed length. 

The authors asked ChatGPT 6 Astra to prove a nonabelian analogue of the random matrix universality theorem of the authors and Yu \cite[Theorem 1.5]{JLY26}. The model was unable to prove the requested analogue. It only proved the weaker nonabelian result in Theorem \ref{thm1d}, which recovers the finite-prime cokernel law of Nguyen and Wood \cite[Theorem 4.1]{NW22} by abelianization.

The authors critically evaluated and revised the models' suggestions, independently verified all mathematical arguments, and wrote the final manuscript. The authors take full responsibility for the correctness of the results.


\end{document}